# CLEAR: Covariant LEAst-square Re-fitting with applications to image restoration

Charles-Alban Deledalle*, Nicolas Papadakis*, Joseph Salmon†, AND Samuel Vaiter‡

**Abstract.** In this paper, we propose a new framework to remove parts of the systematic errors affecting popular restoration algorithms, with a special focus for image processing tasks. Generalizing ideas that emerged for $\ell_1$ regularization, we develop an approach re-fitting the results of standard methods towards the input data. Total variation regularizations and non-local means are special cases of interest. We identify important covariant information that should be preserved by the re-fitting method, and emphasize the importance of preserving the Jacobian (w.r.t. the observed signal) of the original estimator. Then, we provide an approach that has a "twicing" flavor and allows re-fitting the restored signal by adding back a local affine transformation of the residual term. We illustrate the benefits of our method on numerical simulations for image restoration tasks.

**1. Introduction.** Restoring an image from its single noisy and incomplete observation necessarily requires enforcing regularity or a model *prior* on the targeted properties of the sought solution. Regularity properties such as sparsity or gradient sparsity of an image are difficult to enforce in general, and notably lead to combinatorial and non-convex problems. When one is willing to guarantee such kinds of features, convex relaxation is a popular path. This is typically done using the $\ell_1$ norm instead of the $\ell_0$ pseudo-norm, as for the Lasso [44] or the total variation [38]. Nevertheless, such relaxations are well known to create solutions with a larger bias.

Typically, for the Lasso, using the $\ell_1$ convex relaxation of the $\ell_0$ pseudo-norm leads large coefficients to be shrunk towards zero.

For the total variation, the same relaxation on the jump amplitudes induces a loss of contrast in the recovered signal; see Figure 1.(a) for an illustration in this case.

In the Lasso case, a well known re-fitting scheme consists in performing *a posteriori* a least-square re-estimation of the non-zero coefficients of the solution. This post re-fitting technique became popular under various names in the literature: Hybrid Lasso [19], Lasso-Gauss [36], OLS post-Lasso [1], Debiased Lasso (see [28, 1] for extensive details on the subject). For the anisotropic total-variation (aniso-TV), the same post re-fitting approach can be performed to re-estimate the amplitudes of the jumps, provided their locations have been correctly identified.

In this paper, we introduce a generalization of this re-fitting technique that aims at re-enhancing the estimation towards the data without altering the desired properties imposed by the model prior. The underlying reason is that if the user choose an estimator for which theoretical guarantees have been proven (as for the Lasso and TV), re-fitting should preserves these guarantees. To that end, we introduce the *Covariant LEAst-square Re-fitting* (CLEAR). Though this method was originally elaborated with $\ell_1$ analysis problems in mind, it has the ability to generalize to a wider family, while in simple cases such as the Lasso or the aniso-TV, it recovers the classical post re-fitting solution described earlier. For instance, our methodology successfully applies to the Tikhonov regularization [46], the isotropic total-variation (iso-TV), the non-local means [3], the block matching 3D (BM3D) [10] and the Dual Domain Image Denoising (DDID) [26]. In common variational contexts, *e.g.,* $\ell_1 - \ell_2$ analysis [20] (encompassing the Lasso, the group Lasso [29, 54, 32], the aniso- and iso-TV), we show that our re-fitting technique can be performed with a complexity overload of about twice that of the original algorithm. In other cases, we introduce a scheme requiring about three times the complexity of the original algorithm.

*IMB, CNRS, Université de Bordeaux, Bordeaux INP, F-33405 Talence, France ({charles-alban.deledalle, nicolas.papadakis}@math.u-bordeaux.fr).

†LTCI, CNRS, Télécom ParisTech, Université Paris-Saclay, 75013 Paris, France (joseph.salmon@telecom-paristech.fr).

‡IMB, CNRS, Université de Bourgogne, 21078 Dijon, France (samuel.vaiter@u-bourgogne.fr).

While our covariant re-fitting technique recovers the classical post re-fitting solution in specific cases, the proposed algorithm offers more stable solutions. Unlike the Lasso post re-fitting technique, ours does not require identifying *a posteriori* the support of the solution, *i.e.,* the set of non-zero coefficients. In the same vein, it does not require identifying the jump locations of the aniso-TV solution. Since the Lasso or the aniso-TV solutions are usually obtained through iterative algorithms, stopped at a prescribed convergence accuracy, the support or jump numerical identification might be imprecise (all the more for ill-posed problems). Such erroneous support identifications lead to results that strongly deviate from the sought re-fitting. Our covariant re-fitting jointly estimates the re-enhanced solution during the iterations of the original algorithm and, as a by product, produces more robust solutions.

This work follows a preliminary study [14] that attempted to suppress the bias emerging from the choice of the method (*e.g.,* $\ell_1$ relaxation), while leaving unchanged the bias due to the choice of the model (*e.g.,* sparsity). While the approach from [14] – hereafter referred to as *invariant re-fitting* – provides interesting results, it is however limited to a class of restoration algorithms that satisfy restrictive local properties.

In particular, the invariant re-fitting cannot handle iso-TV. In this case, the invariant re-fitting is unsatisfactory as it removes some desired aspects enforced by the prior, such as smoothness, and suffers from a significant increase of variance. A simple illustration of this phenomenon for iso-TV is provided in Figure 2.(d) where artificial oscillations are wrongly amplified near the boundaries.

While the covariant and the invariant re-fitting both correspond to the least-square post re-fitting step in the case of aniso-TV, the two techniques do not match for iso-TV. Indeed, CLEAR outputs a more relevant solution than the one from the invariant re-fitting. Figure 2.(e) shows the benefit of our proposed solution *w.r.t.* the (naive) invariant re-fitting displayed in Figure 2.(d).

It is worth mentioning that the covariant re-fitting is also strongly related to boosting methods re-injecting useful information remaining in the residual (*i.e.,* the map of the point-wise difference between the original signal and its prediction). Such approaches can be traced back to *twicing* [47] and have recently been thoroughly investigated: boosting [4], Bregman iterations and nonlinear inverse scale spaces [33, 5, 52, 34], ideal spectral filtering in the analysis sense [23], SAIF-boosting [30, 43] and SOS-boosting [37] being some of the most popular ones. Most of these methods can be performed iteratively, leading to a difficult choice for the number of steps to consider in practice. Our method has the noticeable advantage that it is by construction a two-step one. Iterating more would not be beneficial (see, Subsection 4.2). Unlike re-fitting, these later approaches aim at improving the overall image quality by authorizing the re-enhanced result to deviate strongly from the original biased solution. In particular, they do not recover the aforementioned post re-fitting technique in the Lasso case. Our objective is not to guarantee the image quality to be improved but rather to generalize the re-fitting approach with the ultimate goal of reducing the bias while preserving the structure and the regularity of the original biased solution.

Interestingly, we have also realized that our scheme presents some similarities with the classical shrinking estimators introduced in [41], especially as presented in [22]. Indeed the step performed by CLEAR, is similar to a shrinkage step with a data-driven residual correction weight (later referred to as $\rho$ in our approach, see Definition 14) when performing shrinkage as in [22, Section 3.1].

Last but not least, it is well known that bias reduction is not always favorable in terms of mean square error (MSE) because of the so-called bias-variance trade-off. It is important to highlight that a re-fitting procedure is expected to re-inject part of the variance, therefore it could lead to an increase of residual noise. Hence, the MSE is not expected to be improved by re-fitting techniques (unlike the aforementioned boosting-like methods that attempt to improve the MSE). We will show in our numerical experiments, that re-fitting is in practice beneficial when the signal of interest fits well the model imposed by the prior. In other scenarios, when the model mismatches the sought signal, the original biased estimator remains favorable in terms of MSE. Re-fitting is nevertheless essential in the latter case for applications where the image intensities have a physical sense and critical decisions are taken from their values.

**2. Background models and notation.** We tackle the problem of estimating an unknown vector $x_0 \in \mathbb{R}^p$ from noisy and incomplete observations

$$y = \Phi x_0 + w \ , \tag{1}$$

where $\Phi$ is a linear operator from $\mathbb{R}^p$ to $\mathbb{R}^n$ and $w \in \mathbb{R}^n$ is the realization of a noisy random vector. This linear model is widely used in statistics and in imagery (*e.g.,* for encoding degradations such as entry-wise masking, convolution, etc.). Typically, the inverse problem associated to (1) is ill-posed, and one should add additional information to recover at least an approximation of $x_0$ from $y$. We denote by $\hat{x} : y \mapsto \hat{x}(y)$ the procedure that provides such an estimation of $x_0$. In this paper, we consider a popular class of estimators relying on a variational formulation that involves a data fidelity term $F(x, y)$ and a regularizing term $G(x)$[1]:

$$\hat{x}(y) \in \underset{x \in \mathbb{R}^p}{\operatorname{argmin}} \ F(x, y) + G(x) \ . \tag{2}$$

Another kind of estimator can be defined as the output of an iterative algorithm $(k, y) \mapsto x^k$, *e.g.,* a numerical solver. For a chosen iteration $k$, we define the final estimator $\hat{x}(y) = x^k$. Such a framework includes for instance proximal splitting methods, as well as discretization of partial differential equations, though we do not investigate this latter road in details.

**2.1. Notation.** For a matrix $M$, $M^+$ is its Moore-Penrose pseudo-inverse. For a (closed convex) set $C$, $\Pi_C$ is the Euclidean projection over $C$ and $\iota_C$ is its indicator function defined by $\iota_C(u) = 0$ if $u \in C$, $+\infty$ otherwise.

For any integer $d \in \mathbb{N}^*$, we write $[d]$ for the set $\{1, \ldots, d\}$. For any subset $\mathcal{I}$ of $[d]$, its complement in $[d]$ is written $\mathcal{I}^c$. For any vector $v$, $v_{\mathcal{I}} \in \mathbb{R}^{|\mathcal{I}|}$ is the sub-vector whose elements are indexed by $\mathcal{I} \subset \mathbb{N}$ and $|\mathcal{I}|$ is its cardinality. For any matrix $M$, $M_{\mathcal{I}}$ is the sub-matrix whose columns are indexed by $\mathcal{I}$. We denote respectively by $\operatorname{Im}[A]$ and $\operatorname{Ker}[A]$ the image space and the kernel space of an operator $A$. For any vector $x \in \mathbb{R}^d$, and $q \in [1, +\infty]$ we denote by $\|x\|_q$ its standard $\ell_q$ norm, and by $\|x\|_0$, the number of non-zero entries of $x$, *i.e.,* $\|x\|_0 = |\{i \in [d] \ : \ x_i \neq 0\}|$.

[1] Often, the solution of (2) is non unique, but for simplicity we only consider a selection of such solutions, and we assume that the selected path $\hat{x} : y \mapsto \hat{x}(y)$ is differentiable almost everywhere.

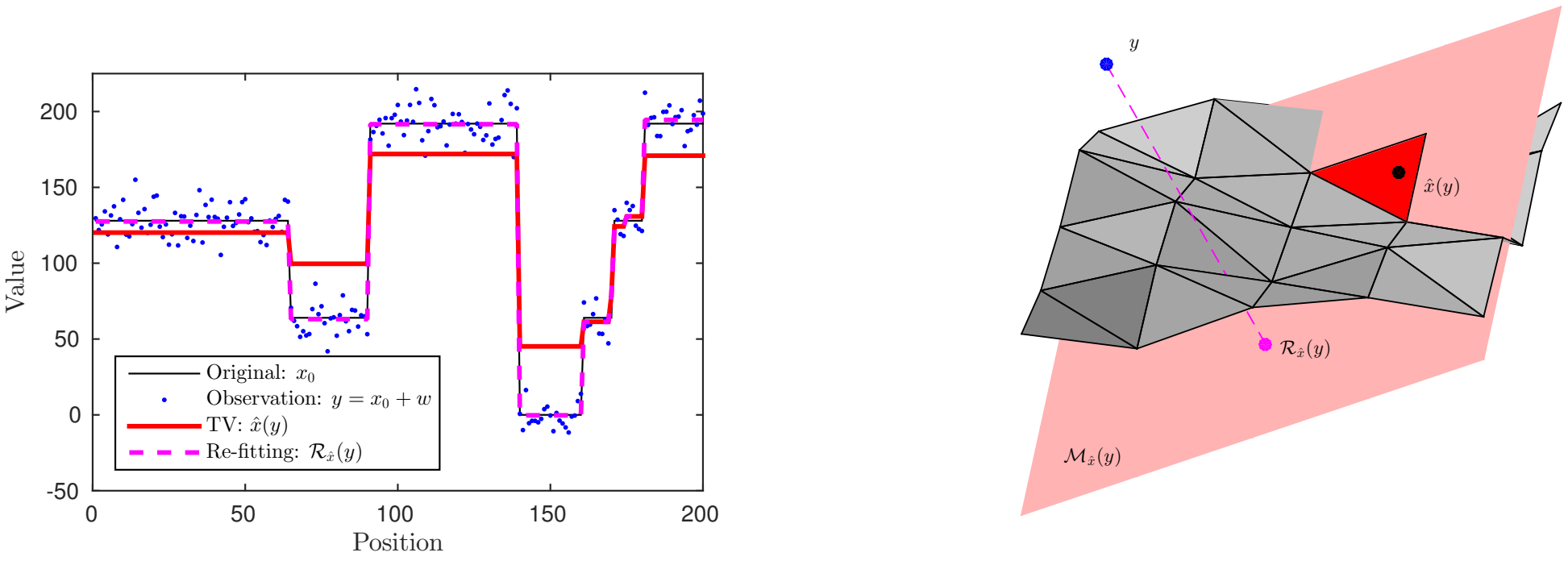

(a) (b)

Fig. 1. *(a) Solutions of 1D-TV and our re-fitting on a noisy signal. (b) Illustration of the invariant re-fitting in a denoising problem of dimension $p = 3$. The gray surface is a piece-wise affine mapping that models the evolution of $\hat{x}$ in an extended neighborhood of $y$. The light red affine plane is the model subspace, i.e., the set of images sharing the same jumps as those of the solution $\hat{x}(y)$. The red triangle is the restriction of the model subspace to images that can be produced by TV. Finally, the pink dot represents the re-fitting $\mathcal{R}^{inv}_{\hat{x}}(y)$ as the orthogonal projection of $y$ on $\mathcal{M}_{\hat{x}}(y)$.*

**2.2. Main investigated estimators.** Here, we provide several canonical examples of estimators of the form (2) that help us illustrating our methodology.

**The affine constrained least-squares**, used when $x_0$ belongs to the affine subspace $C = b + \mathrm{Im}[A]$ with $b \in \mathbb{R}^p$ and $A \in \mathbb{R}^{p \times n}$, are a particular case of (2) where $F(x, y) = \frac{1}{2}\|y - \Phi x\|_2^2$ and $G(x) = \iota_C(x)$. The solution of minimum Euclidean norm is unique and given by: $\hat{x}(y) = b + A(\Phi A)^+(y - \Phi b)$ (see Appendix A.1).

**The Tikhonov regularization** [46] (or Ridge regression [25]), used when $\|\Gamma x_0\|_2$ is small with $\Gamma \in \mathbb{R}^{m \times p}$, is an instance of (2) where $F(x, y) = \frac{1}{2}\|y - \Phi x\|_2^2$ and $G(x) = \frac{\lambda}{2}\|\Gamma x\|_2^2$, for some parameters $\lambda > 0$. Provided $\operatorname{Ker} \Phi \cap \operatorname{Ker} \Gamma = \{0\}$, $\hat{x}(y)$ exists and is uniquely defined as: $\hat{x}(y) = (\Phi^\top \Phi + \lambda \Gamma^\top \Gamma)^{-1} \Phi^\top y$ (see Appendix A.2).

**The hard-thresholding** [16], used when $\Phi = \mathrm{Id}$ and $x_0$ is sparse, is a solution of (2) where $F(x, y) = \frac{1}{2}\|y - x\|_2^2$ and $G(x) = \frac{\lambda^2}{2}\|x\|_0$, for some parameter $\lambda > 0$. The hard-thresholding operation writes: $\hat{x}(y)_{\mathcal{I}} = y_{\mathcal{I}}$ and $\hat{x}(y)_{\mathcal{I}^c} = 0$, where $\mathcal{I} = \{i \in [p] \ : \ |y_i| > \lambda\}$ (see Appendix A.3).

**The soft-thresholding** [16], used when $\Phi = \mathrm{Id}$ and $x_0$ is sparse, is a solution of (2) where $F(x, y) = \frac{1}{2}\|y - x\|_2^2$ and $G(x) = \lambda\|x\|_1$ for some parameter $\lambda > 0$. The soft-thresholding operation writes: $\hat{x}(y)_{\mathcal{I}} = y_{\mathcal{I}} - \lambda \operatorname{sign}(y_{\mathcal{I}})$ and $\hat{x}(y)_{\mathcal{I}^c} = 0$, where $\mathcal{I}$ is defined as for the hard-thresholding (see Appendix A.4).

**The $\ell_1$ synthesis** (or Lasso [44, 16]), used when $x_0$ is sparse, is a solution of (2) where $F(x, y) = \frac{1}{2}\|\Phi y - x\|_2^2$ and $G(x) = \lambda\|x\|_1$ for some parameter $\lambda > 0$. Provided the solution is unique (see for instance [45]), it reads

$$\hat{x}(y)_{\mathcal{I}} = (\Phi_{\mathcal{I}})^+ y_{\mathcal{I}} - \lambda((\Phi_{\mathcal{I}})^\top \Phi_{\mathcal{I}})^{-1} s_{\mathcal{I}} \quad \text{and} \quad \hat{x}(y)_{\mathcal{I}^c} = 0 \ , \tag{3}$$

for almost all $y$ and where $\mathcal{I} = \operatorname{supp}(\hat{x}(y)) = \{i \in [p] \ : \ \hat{x}(y)_i \neq 0\}$ is the support of $\hat{x}(y)$, and $s_{\mathcal{I}} = \operatorname{sign}(\hat{x}(y)_{\mathcal{I}})$.

**The $\ell_1$ analysis**, used when $\Gamma x_0$ is sparse with $\Gamma \in \mathbb{R}^{m \times p}$, is a solution of (2) where

$$F(x, y) = \frac{1}{2}\|\Phi x - y\|_2^2 \quad \text{and} \quad G(x) = \lambda\|\Gamma x\|_1 \ , \quad \text{for some } \lambda > 0 \ . \tag{4}$$

Provided $\operatorname{Ker}\Phi \cap \operatorname{Ker}\Gamma = \{0\}$, there exists a solution given implicitly, see [48], as

$$\hat{x}(y) = U(\Phi U)^{+}y - \lambda U(U^\top \Phi^\top \Phi U)^{-1}U^\top (\Gamma^\top)_{\mathcal{I}} s_{\mathcal{I}} \ , \tag{5}$$

for almost all $y$ and where $\mathcal{I} = \operatorname{supp}(\Gamma\hat{x}(y)) = \{i \in [m] \ : \ (\Gamma\hat{x}(y))_i \neq 0\}$ is called the $\Gamma$-support of the solution, $s_{\mathcal{I}} = \operatorname{sign}((\Gamma\hat{x}(y))_{\mathcal{I}})$, $U$ is a matrix whose columns form a basis of $\operatorname{Ker}[\Gamma_{\mathcal{I}^c}]$ and $\Phi U$ has full column rank.

The anisotropic Total-Variation (aniso-TV) [38] is an instance of (4) where $x_0 \in \mathbb{R}^p$ is identified to a $b$-dimensional signal, for which $\Gamma = \nabla : \mathbb{R}^p \to \mathbb{R}^{p\times b}$ is the discrete gradient operator and $\|\nabla x\|_1 = \sum_{i=1}^p \|(\nabla x)_i\|_1$. Aniso-TV promotes piece-wise constant solutions with large constant regions and few sharp transitions. Here, $\mathcal{I}$ is the set of indexes where the solution has discontinuities. The $\ell_1$ norm of the gradient field induces an anisotropic effect by favoring the jumps to be aligned with the canonical directions (so it favors squared like structures rather than rounded ones).

**The $\ell_1-\ell_2$ analysis** [29, 54, 32], used when $\Gamma x_0$ is block sparse with $\Gamma \in \mathbb{R}^p \mapsto \mathbb{R}^{m\times b}$, is a solution of (2) with

$$F(x,y) = \frac{1}{2}\|\Phi x - y\|_2^2 \quad \text{and} \quad G(x) = \lambda\|\Gamma x\|_{1,2} \ , \tag{6}$$

for some parameter $\lambda > 0$ and where $\|\Gamma x\|_{1,2} = \sum_{i=1}^m \|(\Gamma x)_i\|_2$. The isotropic Total-Variation (iso-TV) [38] is a particular instance of (6) where $x_0 \in \mathbb{R}^p$ can be identified to a $b$-dimensional signal, for which $\Gamma = \nabla : \mathbb{R}^p \to \mathbb{R}^{p\times b}$ and $\|\nabla x\|_{1,2} = \sum_{i=1}^p \|(\nabla x)_i\|_2$. Like aniso-TV, it promotes solution with large constant regions, but some transition regions can be smooth (see, *e.g.,* [6]), typically those with high curvature in the input image $y$, see Figure 2.(a)-(c). A major difference is that the $\ell_1 - \ell_2$ norm induces an isotropic effect by favoring rounded like structures rather than squared ones.

**The (block-wise) non-local means** [3], used when the image $x_0 \in \mathbb{R}^p$ (with $p = p_1 \times p_2$) is composed of many redundant patterns, is the solution of (2) for $F(x,y) = \frac{1}{2}\sum_{i,j} w_{i,j}\|\mathcal{P}_i x - \mathcal{P}_j y\|_2^2$ and $G(x) = 0$, where $\mathcal{P}_i$ is the linear operator extracting a patch (*i.e.,* a small window) at pixel $i$ of size $(2b+1)^2$. The index $i \in [p_1]\times[p_2]$ spans the image domain and $j - i \in [-s,s]^2$ spans a limited search window. The weights $w_{i,j}$ are defined as $w_{i,j} = \varphi\left(\|\mathcal{P}_i y - \mathcal{P}_j y\|_2^2\right)$ where the kernel $\varphi : \mathbb{R}^+ \to [0,1]$ is a decreasing function. Assuming periodical boundary conditons, its solution is given as

$$\hat{x}(y)_i = \frac{\sum_j \bar{w}_{i,j} y_j}{\sum_j \bar{w}_{i,j}} \quad \text{with} \quad \bar{w}_{i,j} = \sum_k w_{i-k,j-k} \ , \tag{7}$$

where $k \in [-d,d]^2$ spans the patch domain (see Appendix A.5).

**2.3. Limitations.** It is important to note that some of the previously introduced estimators are known to suffer from systematic drawbacks.

For instance, the Lasso contracts large coefficients towards 0 by the shift $\lambda((\Phi_{\mathcal{I}})^\top\Phi_{\mathcal{I}})^{-1}s_{\mathcal{I}}$ whose expression simplifies to $\lambda \operatorname{sign}(y_{\mathcal{I}})$ for the soft-thresholding. In the same vein, the $\ell_1$ analysis contracts the signal towards $\operatorname{Ker}[\Gamma]$ by the shift $\lambda U(U^\top\Phi^\top\Phi U)^{-1}U^\top(\Gamma^\top)_{\mathcal{I}}s_{\mathcal{I}}$. In the particular case of aniso-TV, this last quantity is the well known systematic loss of contrast of TV: a shift of intensity on each piece depending on its surrounding and the ratio between its perimeter and its area (see [42] for a thorough study for the 1D case). A simple illustration is provided in Figure 1 where TV is used for denoising a 1D piece-wise constant signal in $[0,192]$ and damaged by additive white Gaussian noise (AWGN) with a standard deviation $\sigma = 10$. Even though TV has perfectly retrieved the support of $\nabla x_0$ with one more extra jump, the intensities of some regions are systematically under- or overestimated. Iso-TV also suffers from a systematic loss of contrast, as illustrated in Figure 2 (c).

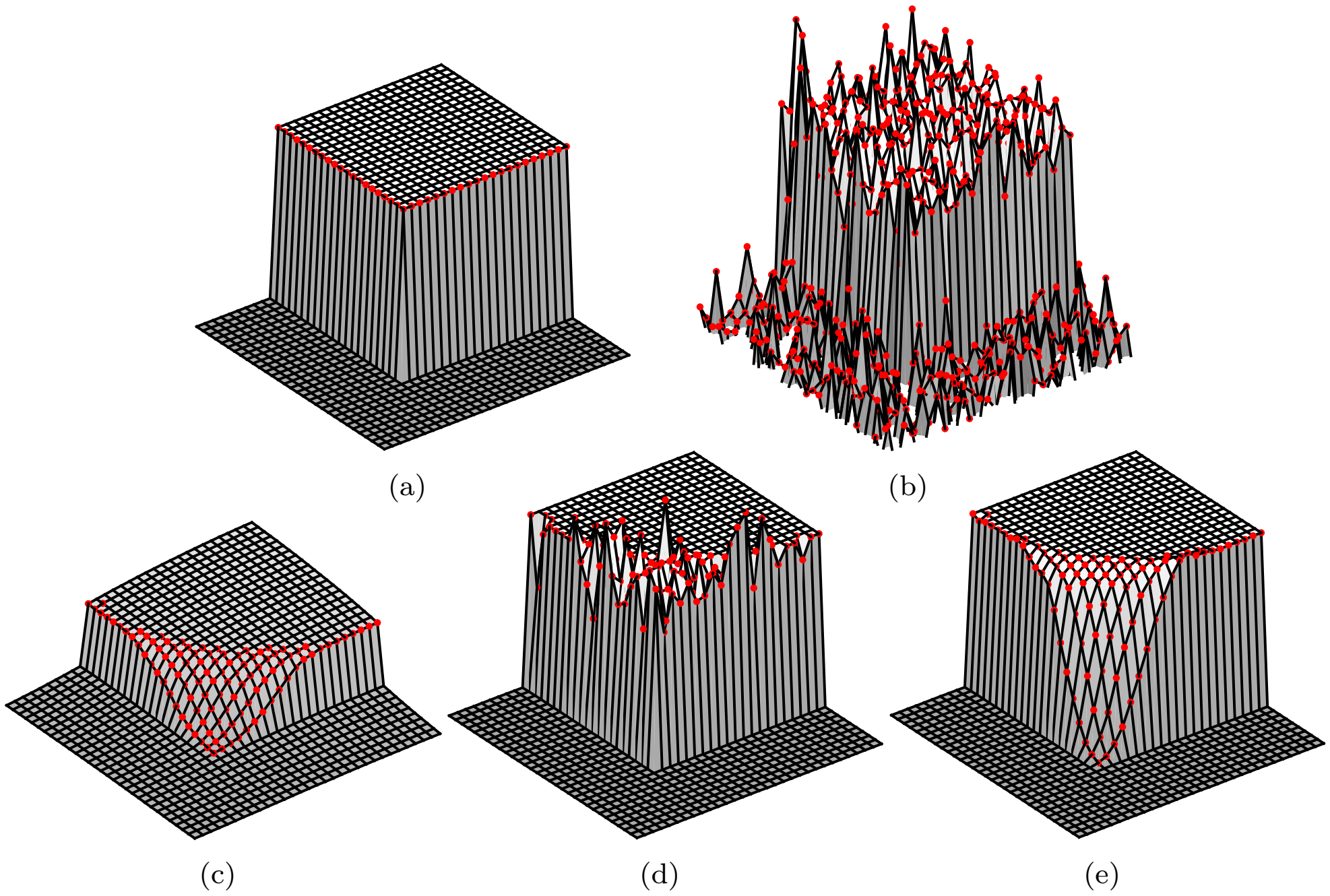

FIG. 2. *(a) A piece-wise constant signal. (b) Its noisy version. (c) Solution of iso-TV with $\lambda = 10$ on the noisy signal. (d) Solution of the invariant re-fitting of iso-TV. (e) Solution of the covariant re-fitting of iso-TV. Red points indicate locations where the discrete gradient is non-zero.*

In the following we investigate possible solutions to reduce such artifacts.

**3. Invariant least-square re-fitting.** As mentioned earlier, practitioners have realized that a systemic contraction affect estimators like the Lasso and the aniso-TV. In the Lasso case, a simple remedy (presented in the introduction) is to perform *a posteriori* a least-square re-fitting step of the non-zero coefficients, *i.e.,* constrained to the support $\mathcal{I}$ of the Lasso solution $\hat{x}(y)$, given by

$$\underset{x\,;\,\mathrm{supp}(x)\subseteq\mathcal{I}}{\mathrm{argmin}} \ \tfrac{1}{2}\|\Phi x - y\|_2^2 \ . \tag{8}$$

In this section, we present a re-fitting procedure, discussed in [14], that generalizes this approach to a broad family of estimators.

**3.1. Re-fitting through model subspace least-squares.** We investigate a re-fitting procedure well suited for estimators almost everywhere differentiable (*a.e.* differentiable). It relies on the notion of model subspace, which requires Jacobian matrix computations. From now on, we consider only estimators $y \mapsto \hat{x}(y)$ from $\mathbb{R}^n$ to $\mathbb{R}^p$ and *a.e.* differentiable.

Many estimation procedures rely on a structural prior of the data. Such structures include smoothness, sparsity, auto-similarity or the fact that the signal is piece-wise constant. Such priors can be captured by the following notion of model subspace.

*Definition* 1. The *model subspace* associated to an *a.e.* differentiable estimator $\hat{x}$ is defined at almost all points $y \in \mathbb{R}^n$ by the affine subspace of $\mathbb{R}^p$

$$\mathcal{M}_{\hat{x}}(y) = \hat{x}(y) + \mathrm{Im}\,[J_{\hat{x}}(y)] \quad , \tag{9}$$

where $J_{\hat{x}}(y) \in \mathbb{R}^{p\times n}$ is the Jacobian matrix of $\hat{x}$ taken at $y$.

The model subspace captures what is linearly invariant through $\hat{x}$ *w.r.t.* small perturbations of $y$, typically, for the Lasso, it will encode the set of signals sharing the same support. In order to generalize the re-fitting step, it is thus natural to cast it as a constrained optimization procedure preserving the model subspace.

*Definition* 2. The *invariant re-fitting* associated to an *a.e.* differentiable estimator $y \mapsto \hat{x}(y)$ is given for almost all $y \in \mathbb{R}^n$ by

$$\mathcal{R}^{\mathrm{inv}}_{\hat{x}}(y) = \hat{x}(y) + J(\Phi J)^{+}(y - \Phi\hat{x}(y)) \in \underset{x\in\mathcal{M}_{\hat{x}}(y)}{\mathrm{argmin}}\ \tfrac{1}{2}\|\Phi x - y\|_2^2 \ , \tag{10}$$

where $J = J_{\hat{x}}(y)$ is the Jacobian matrix of $\hat{x}$ at the point $y$. In the following, we use the notation $J$ when no ambiguity is possible.

*Remark* 3. We only consider $F(x, y) = \frac{1}{2}\|\Phi x - y\|_2^2$, but extensions to a general $F$ (*e.g.,* logistic regression) is straightforward and reads $\mathcal{R}^{\mathrm{inv}}_{\hat{x}}(y) \in \underset{x\in\mathcal{M}_{\hat{x}}(y)}{\mathrm{argmin}}\ F(x, y)$.

*Remark* 4. When $\hat{x}(y) \in \mathrm{Im}[J]$, then $\mathcal{M}_{\hat{x}}(y) = \mathrm{Im}[J]$ and $\mathcal{R}^{\mathrm{inv}}_{\hat{x}}(y) = J(\Phi J)^{+}(y)$.

**3.2. Re-fitting examples.** We now exemplify the previous definitions for the various variational estimators introduced in Subsection 2.2.

**The affine constrained least-squares** have everywhere the same Jacobian matrix $J = A(\Phi A)^{+}$ and its affine model subspace is $\mathcal{M}_{\hat{x}}(y) = b + \mathrm{Im}[A(\Phi A)^{\top}]$ (as $\mathrm{Im}[M^{+}] = \mathrm{Im}[M^{\top}]$). Taking $C = \mathbb{R}^p$, $n = p$, $A = \mathrm{Id}$ and $b = 0$, leads to $\hat{x}(y) = \Phi^{+}y$ whose model subspace is $\mathcal{M}_{\hat{x}}(y) = \mathrm{Im}[\Phi^{\top}]$, reducing to $\mathbb{R}^p$ when $\Phi$ has full column rank. In this case, the invariant re-fitting is $\mathcal{R}^{\mathrm{inv}}_{\hat{x}}(y) = \hat{x}(y)$.

**The Tikhonov regularization** has everywhere the same Jacobian matrix $J = (\Phi^{\top}\Phi + \lambda\Gamma^{\top}\Gamma)^{-1}\Phi^{\top}$ and everywhere the same affine model subspace $\mathcal{M}_{\hat{x}}(y) = \mathrm{Im}[J]$. If follows that $\mathcal{R}^{\mathrm{inv}}_{\hat{x}}(y) = J(\Phi J)^{+}y$. In particular, when $\Phi$ has full column rank, $\mathcal{M}_{\hat{x}}(y) = \mathbb{R}^p$ and $\mathcal{R}^{\mathrm{inv}}_{\hat{x}}(y) = \Phi^{+}y$.

**The soft- and hard-thresholding** share *a.e.* the same Jacobian matrix given by $J = \mathrm{Id}_{\mathcal{I}} \in \mathbb{R}^{p\times|\mathcal{I}|}$. Their model subspace reads as $\mathcal{M}_{\hat{x}}(y) = \mathrm{Im}[\mathrm{Id}_{\mathcal{I}}] = \mathrm{Im}[J]$ and the invariant re-fitting is for both the hard-thresholding.

**The $\ell_1$ synthesis** (or Lasso) has for Jacobian matrix $J = \mathrm{Id}_{\mathcal{I}}(\Phi_{\mathcal{I}})^{+}$ *a.e.* and when $\Phi_{\mathcal{I}}$ has full column rank, it shares the same model subspace $\mathcal{M}_{\hat{x}}(y) = \mathrm{Im}[\mathrm{Id}_{\mathcal{I}}] = \mathrm{Im}[J]$ as the hard- and soft-thresholding. Its invariant re-fitting is in this case $\mathcal{R}^{\mathrm{inv}}_{\hat{x}}(y) = \mathrm{Id}_{\mathcal{I}}(\Phi_{\mathcal{I}})^{+}y$. While the Lasso systematically underestimates the amplitude of the signal by a shift $\lambda\mathrm{Id}_{\mathcal{I}}((\Phi_{\mathcal{I}})^{t}\Phi_{\mathcal{I}})^{-1}s_{\mathcal{I}}$, the re-fitting $\mathcal{R}^{\mathrm{inv}}_{\hat{x}}(y)$ is free of such a contraction.

**The $\ell_1$ analysis** has for Jacobian matrix $J = U(\Phi U)^{+}$ *a.e.* (since the $\Gamma$-support, $\mathcal{I}$ and the sign $s_{\mathcal{I}}$ are *a.e.* stable *w.r.t.* small perturbations [48]). It results that the model subspace reads as $\mathcal{M}_{\hat{x}}(y) = \mathrm{Ker}[\Gamma_{\mathcal{I}^c}] = \mathrm{Im}[J]$. The generalization of the Lasso re-fitting leads to the least-square estimator constrained on the $\Gamma$-support of $\hat{x}(y)$:

$$\mathcal{R}^{\mathrm{inv}}_{\hat{x}}(y) = U(\Phi U)^{+}y \ . \tag{11}$$

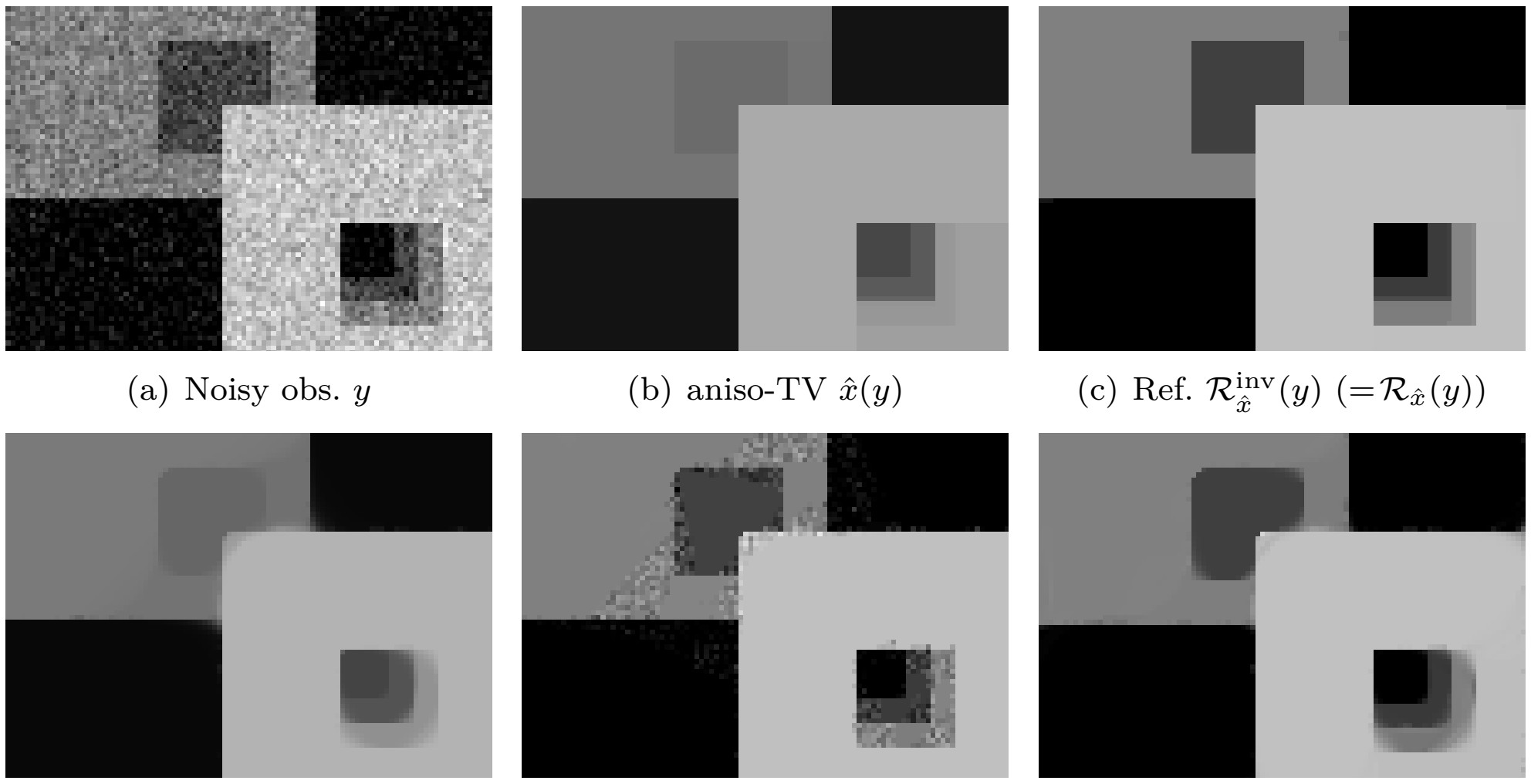

(a) Noisy obs. $y$ (b) aniso-TV $\hat{x}(y)$ (c) Ref. $\mathcal{R}^{\text{inv}}_{\hat{x}}(y)$ $(=\mathcal{R}_{\hat{x}}(y))$

(d) iso-TV $\hat{x}(y)$ (e) Inv. ref. $\mathcal{R}^{\text{inv}}_{\hat{x}}(y)$ (f) Cov. ref. $\mathcal{R}_{\hat{x}}(y)$

FIG. 3. *(a) A noisy image $y = x_0 + w$ where $x_0$ contains six overlapping squares. (b) The solution $\hat{x}(y)$ of aniso-TV with $\lambda = 10$, and(c) its invariant re-fitting $\mathcal{R}^{inv}_{\hat{x}}(y)$ (which coincides in this case with the covariant one $\mathcal{R}_{\hat{x}}(y)$). (d) The solution $\hat{x}(y)$ of iso-TV with $\lambda = 10$, (e) its invariant re-fitting $\mathcal{R}^{inv}_{\hat{x}}(y)$ and (f) its covariant one $\mathcal{R}_{\hat{x}}(y)$.*

For aniso-TV denoising (*i.e.,* $\Phi = \text{Id}$ and $\Gamma = \nabla$), the model subspace is the space of images whose jumps are included in those of the solution $\hat{x}(y)$. The re-fitting procedure $\mathcal{R}^{\text{inv}}_{\hat{x}}$ is the projector $\Pi_{\text{Im}[J]} = UU^+$ that performs a piece-wise average of its input on each plateau of the solution (see Figure 1.(a) for an illustration in 1D).

**The $\ell_1 - \ell_2$ analysis** admits *a.e.* the following Jacobian operator [49]

$$(12) \qquad J = U(U^\top \Phi^\top \Phi U + \lambda U^\top \Gamma^\top \Omega \Gamma U)^{-1} U^\top \Phi^\top y$$

$$\text{where} \quad \Omega : z \in \mathbb{R}^{m \times b} \mapsto \begin{cases} \frac{1}{\|(\Gamma\hat{x}(y))_i\|_2} \left( z_i - \left\langle z_i, \frac{(\Gamma\hat{x}(y))_i}{\|(\Gamma\hat{x}(y))_i\|_2} \right\rangle \frac{(\Gamma\hat{x}(y))_i}{\|(\Gamma\hat{x}(y))_i\|_2} \right) & \text{if} \quad i \in \mathcal{I}, \\ 0 & \text{otherwise,} \end{cases}$$

with $U$ and $\mathcal{I}$ defined exactly as for the $\ell_1$ analysis case. Note that Eq. (12) is well founded as soon as $\text{Ker}[\Gamma_{\mathcal{I}^c}] \cap \text{Ker}\,\Phi = \{0\}$. For weaker conditions, see [49, example 26]. In the particular case where $\Phi U$ has full column rank, the model subspace matches with the one of the $\ell_1$ analysis, *i.e.,* $\mathcal{M}_{\hat{x}}(y) = \text{Ker}[\Gamma_{\mathcal{I}^c}] = \text{Im}[J]$, and the re-fitting is also given by Eq. (11), hence $\mathcal{R}^{\text{inv}}_{\hat{x}}(y) = U(\Phi U)^+ y$. As a consequence, for the iso-TV denoising, *i.e.,* when $\Phi = \text{Id}$ and $\Gamma = \nabla$, the invariant re-fitting step consists again in performing a piece-wise average on each plateau of the solution.

**The non-local means** has a Jacobian matrix with a complex structure [50, 18]. In particular, computing the projection on the model subspace is challenging in this case and so is the computation of the invariant re-fitting. Note that a greedy procedure was proposed in [14] to compute the invariant re-fitting.

**3.3. Results and limitations.** Figure 1.(a) illustrates the invariant re-fitting in the case of a 1D total-variation denoising example ($\ell_1$ analysis estimator). It recovers the jumps of the underlying signal (adding an extra one), but systematically under-estimates their amplitudes. As expected, re-fitting re-enhances the amplitudes

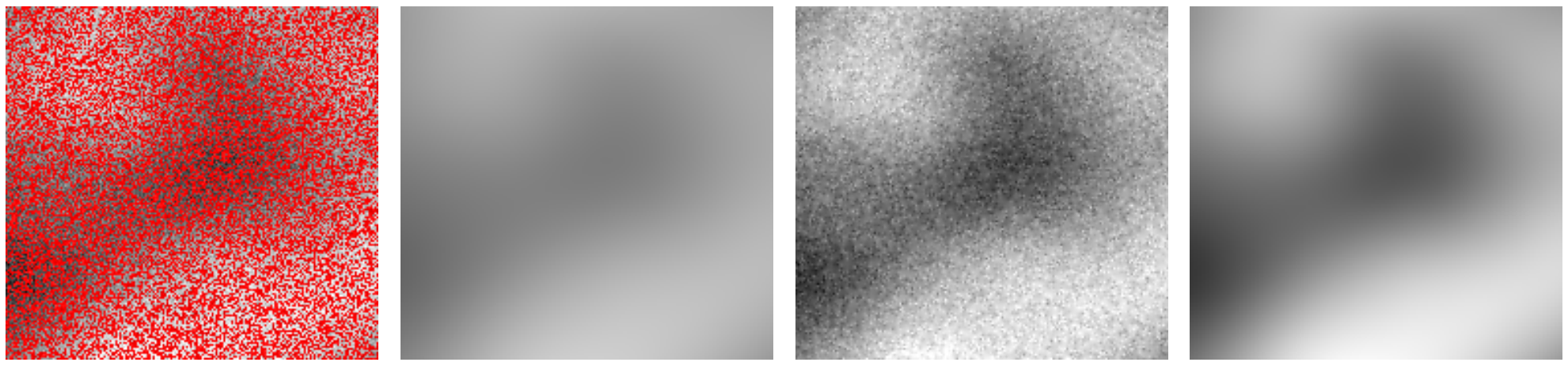

(a) Masked obs. $y$ (b) Tikhonov $\hat{x}(y)$ (c) Inv. ref. $\mathcal{R}^{\text{inv}}_{\hat{x}}(y)$ (d) Cov. ref. $\mathcal{R}_{\hat{x}}(y)$

Fig. 4. *(a) A noisy and incomplete image $y = \Phi x_0 + w$ where $\Phi$ is a masking operator encoding missing pixel values (in red) and $x_0$ is a smooth signal. (b) Tikhonov regularization $\hat{x}(y)$ with $\Gamma = \nabla$ and $\lambda = 20$, (c) its invariant re-fitting $\mathcal{R}^{inv}_{\hat{x}}(y)$ and (d) its covariant one $\mathcal{R}_{\hat{x}}(y)$.*

of all plateaus towards the data. Figure 1.(b) gives a geometrical interpretation in dimension $p = 3$ of the model subspace and the invariant re-fitting. The model subspace is represented as the tangent plane of $\hat{x}$ at $y$ and its re-fitting is the projection of $y$ on this plane. All elements of this plane share the same jumps with the solution $\hat{x}(y)$. Figure 3.(a)-(c) gives a similar illustration for a 2D aniso-TV denoising example.

While the invariant re-fitting acts properly for the $\ell_1$ analysis estimator, it is however less appealing in other scenarios. Figure 2.(c),(d) and Figure 3.(d),(e) give two illustrations of the invariant re-fitting of a 2D iso-TV denoising example. As for aniso-TV, the invariant re-fitting is the projection of $y$ on the space of signals whose jumps are located at the same position as those of $\hat{x}(y)$. But unlike the anisotropic case, $\hat{x}(y)$ is not piece-wise constant. Instead of being composed of distinct flat regions, it reveals smoothed transitions with dense supports (referred to as extended supports in [6]), see Figure 2.(c). As a consequence, the invariant re-fitting re-introduces a large amount of the original noisy signal in these smooth but non-constant areas, creating the artifacts observed on Figure 2.(d) and Figure 3.(e).

Figure 4,(a)-(c) gives another illustration of the invariant re-fitting of a 2D Tikhonov masking example (with $\Phi$ a diagonal matrix with 0 or 1 elements on the diagonal and $\Gamma = \nabla$). While the dynamic of the Tikhonov solution $\hat{x}(y)$ has been strongly reduced, the re-fitting $\mathcal{R}^{\text{inv}}_{\hat{x}}(y)$ re-fits clearly the solution towards the original intensities. However, such a re-fitting is not satisfactory as it does not preserve the smoothness of the solution $\hat{x}(y)$.

In fact, the model subspace captures only what is linearly invariant through $\hat{x}$ *w.r.t.* small perturbations of $y$. This includes the support of the solution for the iso-TV, and the absence of variations inside $\text{Im}[J]^\perp$ for the Tikhonov regularization. In particular, it fails at capturing some of the desirable relationships between the entries of $y$ and the entries of $\hat{x}(y)$, what we call the *covariants*. These relationships typically encode some of the local smoothness and non-local interactions between the entries of the solution $\hat{x}(y)$. Such crucial information is not encoded in the linear model subspace, but interestingly the Jacobian matrix captures by definition how much the entries of $\hat{x}$ linearly varies *w.r.t.* all the entries of $y$. This is at the heart of *the covariant re-fitting* defined in the next section and, for comparison, it produces the results given in Figure 2.(e), Figure 3.(f) and Figure 4.(d).

**4. Covariant LEast-square Re-fitting (CLEAR).** The objective of this section is to present our main contribution, the introduction of the covariant re-fitting procedure (CLEAR). We particularly aim at solving the issues raised in Subsection 3.3. Toward this goal, we put a stronger emphasis on the first-order behavior of the original estimator by imposing the conservation of its Jacobian, at least locally. To define the re-fitting procedure we first need to introduce a procedure which, loosely speaking, takes as input the original estimator and a guess of $\Phi x_0$ and outputs a new estimator with some targeted properties. We next define our Covariant LEast-square re-fitting (CLEAR) by choosing a specific guess for $\Phi x_0$, see Definition 14.

**4.1. Local approach and desired properties for a suitable re-fitting.** In this subsection, our objective is to define, from the original estimator $\hat{x}$ and a guess $z \in \mathbb{R}^n$ of $\Phi x_0$, a new estimator $\mathcal{D}_{\hat{x},z} : \mathbb{R}^n \to \mathbb{R}^p$ that satisfies several desirable properties and shares with $\hat{x}$ some first-order properties. After-wise, we will consider the choice $z = y$, and the resulting estimator is going to be our covariant re-fitting $\mathcal{R}_{\hat{x}}$. We are now equipped to introduce such a guess based re-fitting.

*Definition* 5. Let $\hat{x} : \mathbb{R}^n \to \mathbb{R}^p$ be differentiable at $z \in \mathbb{R}^n$. An estimator $\mathcal{D}_{\hat{x},z}$ : $\mathbb{R}^n \to \mathbb{R}^p$ is a guess based covariant least-square re-fitting of $\hat{x}$ for $z$, if

$$\mathcal{D}_{\hat{x},z} \in \underset{h \in \mathcal{H}}{\operatorname{argmin}} \|\Phi h(z) - z\|_2^2 \ , \tag{13}$$

where $\mathcal{H}$ is the set of maps $h : \mathbb{R}^n \to \mathbb{R}^p$ satisfying, for all $y \in \mathbb{R}^n$,

1. **Affine map:** $h(y) = Ay + b$ for some $A \in \mathbb{R}^{p\times n}, b \in \mathbb{R}^p$,
2. **Covariant preserving:** $J_h(z) = \rho J_{\hat{x}}(z)$ for some $\rho \in \mathbb{R}$,
3. **Coherent map:** $h(\Phi\hat{x}(z)) = \hat{x}(z)$.

Definition 5 is natural as it states that a guess based re-fitting of $\hat{x}$ for $z$ should be, in prediction, as close as possible to $z$. Of course, it should satisfy some extra conditions. First, the estimator should be easy to compute, and so we choose a first order approximation, leading to a locally affine estimator. Second, as highlighted in Subsection 3.3, the relative variation of the original estimator *w.r.t.* the input should be preserved to capture, not only the invariant features of the estimator but also its first-order behavior, capturing both its singularities and smoothness. Third, applying a re-fitting step to the prediction obtained by the original estimator at $z$ should not modify it. The purpose of re-fitting is to be close to $y$, while also preserving the structure of $\hat{x}(z)$. Hence, if $y = \Phi\hat{x}(z)$, the result should be unaltered.

The next theorem provides a unique closed form expression for $\mathcal{D}_{\hat{x},z}(y)$.

THEOREM 6. *Let $\hat{x}$ be an estimator from $\mathbb{R}^n$ to $\mathbb{R}^p$ differentiable at $z \in \mathbb{R}^p$. Then, for $\delta = z - \Phi\hat{x}(z)$, the guess based covariant least-square re-fitting, defined in Definition 5, exists, is unique if $\Phi J\delta \neq 0$, and is given by*

$$\mathcal{D}_{\hat{x},z}(y) = \hat{x}(z) + \rho J(y - \Phi\hat{x}(z)) \quad \textit{where} \quad \rho = \begin{cases} \dfrac{\langle \Phi J\delta, \, \delta\rangle}{\|\Phi J\delta\|_2^2} & \textit{if} \quad \Phi J\delta \neq 0 \ , \\ 1 & \textit{otherwise} \ , \end{cases} \tag{14}$$

*where $J = J_{\hat{x}}(z)$ is the Jacobian matrix of $\hat{x}$ at the point $z$.*

*Proof.* Let $h$ be a mapping satisfying properties 1., 2. and 3. in the previous definition. Observe that properties 1. and 2. of the set $\mathcal{H}$ together ensures that the estimator is of the form $h(y) = \rho Jy + b$ for some $\rho \in \mathbb{R}$ and $b \in \mathbb{R}^p$. Plugging condition 3. gives that $b = (\mathrm{Id} - \rho J\Phi)\hat{x}(z)$, hence $h(y) = \hat{x}(z) + \rho J(y - \Phi\hat{x}(z))$. Reciprocally, it is easy to see that any estimator of the form $h(y) = \hat{x}(z) + \rho J(y - \Phi\hat{x}(z))$ satisfies properties 1., 2. and 3. It thus remains to find $\rho$.

Remark that Problem (13) can be recast as a one-dimensional problem

$$\min_{\rho\in\mathbb{R}} \left\{ \quad \|\Phi(\hat{x}(z)+\rho J(z-\Phi\hat{x}(z)))-z\|_2^2 = \|(\mathrm{Id}-\rho\Phi J)(\Phi\hat{x}(z)-z)\|_2^2 \quad \right\} \tag{15}$$

whose unique solution, if $\Phi J(\Phi\hat{x}(z)-z)\neq 0$ is given by (14) and $\rho=1$ otherwise. □

The case where $\Phi J$ is an orthogonal projector leads to interesting properties for instance when $\hat{x}$ is associated to the constrained least-squares, the Lasso or aniso-TV.

Proposition 7. *If $\Phi J$ is an orthogonal projector, then $\rho=1$.*

*Proof.* If $\Phi J$ is an orthogonal projector, it follows that $\|\Phi J\delta\|_2^2 = \langle \Phi J\delta,\, \Phi J\delta\rangle = \langle (\Phi J)^\top \Phi J\delta,\, \delta\rangle = \langle \Phi J\delta,\, \delta\rangle$. Injecting this in (14) gives $\rho=1$. □

**Statistical interpretation.** For a random vector with expectation $\Phi x_0$ and finite second order moment, bias and covariance of $\mathcal{D}_{\hat{x},z}$ are given in closed form.

Proposition 8. *Let $Y$ be a random vector in $\mathbb{R}^n$ such that $\mathbb{E}[Y]=\Phi x_0$, and $\mathrm{Cov}[Y]=\Sigma\in\mathbb{R}^{n\times n}$. Then $y\mapsto\mathcal{D}_{\hat{x},z}(y)$ satisfies*

$$\mathbb{E}[\mathcal{D}_{\hat{x},z}(Y)]-x_0 = (\mathrm{Id}-\rho J\Phi)(\hat{x}(z)-x_0)\ , \tag{16}$$

$$\mathrm{Cov}[\mathcal{D}_{\hat{x},z}(Y),Y] = \rho J\Sigma\ , \tag{17}$$

$$\mathrm{Cov}[\mathcal{D}_{\hat{x},z}(Y)] = \rho^2 J\Sigma J^\top, \tag{18}$$

*where the cross covariance is $\mathrm{Cov}[X,Y]=\mathbb{E}[XY^\top]-\mathbb{E}[X]\mathbb{E}[Y]^\top$, for any random column vectors $X$ and $Y$ (not necessarily of the same size), and $\mathrm{Cov}[Y]=\mathrm{Cov}[Y,Y]$.*

*Proof.* The first equality is a direct consequence of the linearity of the expectation operator. The second equality arises from the following

$$\begin{aligned}\mathbb{E}[(\hat{x}(z)+\rho J(Y-\Phi\hat{x}(z)))Y^\top]-&\mathbb{E}[(\hat{x}(z)+\rho J(Y-\Phi\hat{x}(z)))]\mathbb{E}[Y]^\top\\ &= \rho J\left(\mathbb{E}[YY^\top]-\mathbb{E}[Y]\mathbb{E}[Y]^\top\right)\end{aligned} \tag{19}$$

since $J$ and $\rho$ are constant *w.r.t.* $y$ as they depend only on the guess $z$. The third equation follows the same sketch by expanding the expression of $\mathrm{Cov}[\mathcal{D}_{\hat{x},z}(Y)]$. □

Proposition 8 provides a closed form expression for the bias, the cross-covariance and the covariance of $\mathcal{D}_{\hat{x},z}$. These quantities are much more intricate to derive for a non-linear estimator $\hat{x}$. Nevertheless, the next corollary shows how these quantities relate to those of the first order Taylor expansion of the original estimator $\hat{x}$.

Corollary 9. *Let $\mathcal{T}_{\hat{x},z}(y)$ be the tangent estimator of $\hat{x}$ at $z\in\mathbb{R}^n$ defined as*

$$\mathcal{T}_{\hat{x},z}(y) = \hat{x}(z)+J(y-z)\ . \tag{20}$$

*Let $Y$ be a random vector in $\mathbb{R}^n$ such that $\mathbb{E}[Y]=\Phi x_0$ and $\mathrm{Cov}[Y]=\Sigma$. Then $y\mapsto\mathcal{T}_{\hat{x},z}(y)$ and $y\mapsto\mathcal{D}_{\hat{x},z}(y)$ satisfy*

$$\mathbb{E}[\mathcal{T}_{\hat{x},z}(Y)]-x_0 = (\hat{x}(z)-x_0)+J(\Phi x_0-z)\ , \tag{21}$$

$$\mathrm{Cov}[\mathcal{D}_{\hat{x},z}(Y),Y] = \rho\mathrm{Cov}[\mathcal{T}_{\hat{x},z}(Y),Y]\ , \tag{22}$$

$$\mathrm{Cov}[\mathcal{D}_{\hat{x},z}(Y)] = \rho^2\mathrm{Cov}[\mathcal{T}_{\hat{x},z}(Y)]\ , \tag{23}$$

$$\mathrm{Corr}[\mathcal{D}_{\hat{x},z}(Y),Y] = \mathrm{Corr}[\mathcal{T}_{\hat{x},z}(Y),Y]\ , \tag{24}$$

$$\mathrm{Corr}[\mathcal{D}_{\hat{x},z}(Y)] = \mathrm{Corr}[\mathcal{T}_{\hat{x},z}(Y)]\ , \tag{25}$$

*where the cross correlation matrix is defined as* $\mathrm{Corr}[X,Y]_{i,j} = \mathrm{Cov}[X,Y]_{i,j}/\sqrt{\mathrm{Cov}[X]_{i,i}\mathrm{Cov}[Y]_{j,j}}$*, for any random column vectors* $X$ *and* $Y$ *(not necessarily of the same size), and* $\mathrm{Corr}[Y] = \mathrm{Corr}[Y,Y]$.

*Proof.* The first relation holds from the expression of $\mathcal{T}_{\hat{x},z}$ and that $J$ does not depend on $y$. It follows that $\mathrm{Cov}[\mathcal{T}_{\hat{x},z}(Y), Y] = J\Sigma$ and $\mathrm{Cov}[\mathcal{T}_{\hat{x},z}(Y)] = J\Sigma J^\top$. These, jointly with Proposition 8, conclude the proof. □

Corollary 9 is essential in this work as it states that, by preserving the Jacobian structure, $\mathcal{D}_{\hat{x},z}(Y)$ cannot depart from the tangent estimator of $\hat{x}$ at $z$ in terms of (cross-)correlations. As a consequence, one can expect that they only differ in terms of expectation, *i.e.*, in terms of bias. The next propositions state that when $\Phi J$ is a projector, the bias in prediction is guaranteed to be reduced by our re-fitting.

Proposition 10. *Let* $Y$ *be a random vector of* $\mathbb{R}^n$ *such that* $\mathbb{E}[Y] = \Phi x_0$*. Assume* $\Phi J$ *is an orthogonal projector, then* $y \mapsto \mathcal{D}_{\hat{x},z}(y)$ *satisfies*

$$\|\Phi(\mathbb{E}[\mathcal{D}_{\hat{x},z}(Y)] - x_0)\|_2 \leqslant \|\Phi(\mathbb{E}[\mathcal{T}_{\hat{x},z}(Y)] - x_0)\|_2 \ . \tag{26}$$

*Proof.* As $\Phi J$ is an orthogonal projector, by virtue of Proposition 7, $\rho = 1$, then

$$\begin{aligned} \|\Phi(\mathbb{E}[\mathcal{T}_{\hat{x},z}(Y)] - x_0)\|_2^2 &= \|\Phi(\hat{x}(z) + J(\Phi x_0 - z) - x_0)\|_2^2 \\ &= \|\Phi(\hat{x}(z) + J\Phi(x_0 - \hat{x}(z)) - x_0)\|_2^2 + \|\Phi J(\Phi\hat{x}(z) - z)\|_2^2 \\ &= \|\Phi(\mathbb{E}[\mathcal{D}_{\hat{x},z}(Y)] - x_0)\|_2^2 + \|\Phi J\delta\|_2^2 \end{aligned} \tag{27}$$

which concludes the proof. □

Proposition 10 is a bit restrictive as it requires $\Phi J$ to be a projector. Nevertheless, this assumption can be relaxed when $z$ satisfies a more technical assumption as shown in the next proposition.

Proposition 11. *Let* $Y$ *be a random vector of* $\mathbb{R}^n$ *such that* $\mathbb{E}[Y] = \Phi x_0$*. Let* $\rho_0 = \frac{\langle \delta_0, \Phi J \delta_0\rangle}{\|\Phi J\delta_0\|_2^2}$ *and* $\delta_0 = \Phi(x_0 - \hat{x}(z))$*. Assume there exists* $\alpha \in [0,1]$ *such that*

$$\left|\frac{\rho - \rho_0}{\rho_0}\right| \leqslant \sqrt{1-\alpha} \ , \tag{28}$$

$$\textit{and} \quad \|\Phi J(\delta - \delta_0)\|_2^2 + 2\langle \delta_0,\, \Phi J(\delta - \delta_0)\rangle \geqslant -\alpha \frac{\langle \delta_0,\, \Phi J\delta_0\rangle^2}{\|\Phi J\delta_0\|_2^2} \ . \tag{29}$$

*Then,* $y \mapsto \mathcal{D}_{\hat{x},z}(y)$ *satisfies*

$$\|\Phi(\mathbb{E}[\mathcal{D}_{\hat{x},z}(Y)] - x_0)\|_2 \leqslant \|\Phi(\mathbb{E}[\mathcal{T}_{\hat{x},z}(Y)] - x_0)\|_2 \ . \tag{30}$$

*Proof.* It follows from Proposition 8 and Corollary 9 that $\|\Phi(\mathbb{E}[\mathcal{D}_{\hat{x},z}(Y)] - x_0)\|_2 = \|(\mathrm{Id} - \rho\Phi J)\delta_0\|_2$ and $\|\Phi(\mathbb{E}[\mathcal{T}_{\hat{x},z}(Y)] - x_0)\|_2 = \|\delta_0 + \Phi J(\delta - \delta_0)\|_2$. Subsequently, we get that Equation (30) holds true if

$$\|(\mathrm{Id} - \rho\Phi J)\delta_0\|_2^2 \leqslant \|\delta_0 + \Phi J(\delta - \delta_0)\|_2^2 \ , \tag{31}$$

$$\textit{i.e.,} \quad \rho^2\|\Phi J\delta_0\|_2^2 - 2\rho\langle \delta_0,\, \Phi J\delta_0\rangle \leqslant \|\Phi J(\delta - \delta_0)\|_2^2 + 2\langle \delta_0,\, \Phi J(\delta - \delta_0)\rangle. \tag{32}$$

Using Assumption (29), a sufficient condition for Equation (30) to hold is

$$\rho^2\|\Phi J\delta_0\|_2^2 - 2\rho\langle \delta_0,\, \Phi J\delta_0\rangle_2 + \alpha\frac{\langle \delta_0,\, \Phi J\delta_0\rangle^2}{\|\Phi J\delta_0\|_2^2} \leqslant 0 \ . \tag{33}$$

The roots of this polynomial are given by $(1 \pm \sqrt{1-\alpha})\rho_0$, which concludes the proof. □

*Remark* 12. Remark that requiring (28) is quite natural as it states that $\rho$ should be close enough to the optimal $\rho_0$ minimizing the discrepancy with regards to $\Phi x_0$ (*i.e.,* minimizing $\|\Phi h(z) - \Phi x_0\|_2^2$ for $h \in \mathcal{H}$ defined as in Definition 5). While the condition (29) sounds more technical, it however holds true in several interesting cases. For instance, when $z = \Phi x_0$, Assumption (29) holds true as it would read $0 \geqslant -\alpha \langle \delta_0, \Phi J \delta_0 \rangle^2 / \|\Phi J \delta_0\|_2^2$ (since $\delta = \delta_0$ and $\rho = \rho_0$). Another case of interest is when $\Phi J$ is an orthogonal projector for which (29) holds true as it would read $\|\Phi J \delta\|_2^2 - \|\Phi J \delta_0\|_2^2 \geqslant -\|\Phi J \delta_0\|_2^2$ (using that $\rho = \rho_0 = 1$, $\langle \cdot, \Phi J \cdot \rangle = \|\Phi J \cdot\|_2^2$, and choosing $\alpha = 1$). Hence, Proposition 11 recovers Proposition 10.

*Remark* 13. Using the same sketch of proof as Proposition 11, the condition $|(\rho - \rho_0)/\rho_0| \leqslant 1$ is sufficient to get $\|\Phi(\mathbb{E}[\mathcal{D}_{\hat{x},z}(Y)] - x_0)\|_2 \leqslant \|\Phi(\hat{x}(z) - x_0)\|_2$. In other words, even though $\rho$ has a relative error of 100% *w.r.t.* $\rho_0$, the estimator $y \mapsto \mathcal{D}_{\hat{x},z}(y)$ still reduces the bias of the constant estimator $y \mapsto \hat{x}(z)$. This result remains valid when comparing $y \mapsto \mathcal{D}_{\hat{x},z}(y)$ to the pseudo-oracle estimator $y \mapsto \hat{x}(z) + J(y - \Phi x_0)$, with the notable difference that they moreover share the same correlation structure.

While it is difficult to state a general result, we can reasonably claim from Proposition 10, Proposition 11 and Remark 12 that the bias tends to be reduced by our re-fitting providing $\Phi J$ is almost a projector (*i.e.,* has eigenvalues concentrated around 0 and 1) and/or $z$ is not too far from $\Phi x_0$. In such cases, the estimator $\mathcal{D}_{\hat{x},z}$ can be considered as a debiasing procedure of $\hat{x}$, in the sense that it reduces the bias of $\mathcal{T}_{\hat{x},z}$ while preserving its correlation structure (according to Corollary 9).

**4.2. Definitions and properties.** Using $\mathcal{D}_{\hat{x},z}$ defined in Theorem 6, we can now give an explicit definition of CLEAR as $\mathcal{R}_{\hat{x}}(y) = \mathcal{D}_{\hat{x},y}(y)$.

*Definition* 14 (CLEAR). The *Covariant LEast-square Re-fitting* associated to an *a.e.* differentiable estimator $y \mapsto \hat{x}(y)$ is, for almost all $y \in \mathbb{R}^n$, given by

$$\mathcal{R}_{\hat{x}}(y) = \hat{x}(y) + \rho J(y - \Phi \hat{x}(y)) \quad \text{with} \quad \rho = \begin{cases} \dfrac{\langle \Phi J \delta, \delta \rangle}{\|\Phi J \delta\|_2^2} & \text{if} \quad \Phi J \delta \neq 0 \ , \\ 1 & \text{otherwise} \ , \end{cases} \tag{34}$$

where $\delta = y - \Phi \hat{x}(y)$ and $J = J_{\hat{x}}(y)$ is the Jacobian matrix of $\hat{x}$ at the point $y$.

Figure 5 gives a geometrical interpretation of CLEAR for a denoising task in dimension $p = 3$. One can observe that if $Y$ varies isotropically, so will its projection on the model subspace. Contrarily, the tangent estimator at a guess $z$ can present an anisotropic behavior along the model subspace, and the guess based re-fitting, which is closer to $z$, will respect this anisotropy in order to capture the local regularity of $\hat{x}$. Finally, the covariant re-fitting is obtained from the guess based re-fitting at $z = y$. For clarity, we assumed that $\mathcal{M}_{\hat{x}}(z) = \mathcal{M}_{\hat{x}}(y)$ on this illustration.

*Remark* 15. The covariant re-fitting performs an additive correction of $\hat{x}(y)$ with a fraction of the directional derivative $J\delta$ in the direction of the residual $\delta$.

*Remark* 16. Observe that in Definition 14, the value of $\rho$ varies when $y$ varies, contrary to the map $y \mapsto \mathcal{D}_{\hat{x},z}(y)$ for which $\rho$ is constant. Note that the mapping $y \mapsto \mathcal{D}_{\hat{x},z}(y)$ is affine, but *not* the map $y \mapsto \mathcal{R}_{\hat{x}}(y)$. Note that, as a consequence, the statistical interpretations given in the previous section do not hold for $\mathcal{R}_{\hat{x}}(y)$ even though they shed some light on its behavior.

**Why not iterating the procedure?** The re-fitting procedure computes $\tilde{x}^2 = \tilde{x}^1 + \rho J(y - \Phi \tilde{x}^1)$, with $\tilde{x}^1 = \hat{x}(y)$. One may wonder if it is beneficial to iterate the process as $\tilde{x}^{k+1} = \tilde{x}^k + \rho J(y - \Phi\tilde{x}^k)$ (in the same vein as [47, 4, 33, 43, 37]). Consider a denoising problem $\Phi = \mathrm{Id}$ with Tikhonov or iso-TV, for which $J$ is symmetrical and $\hat{x}(y) \in \mathrm{Im}[J]$ (see examples of Subsection 3.2). The sequence will converge if and only if $J(y - \tilde{x}^k)$ vanishes, *i.e.,* $\tilde{x}^k$ must converge to $J^+Jy + \zeta$ with $\zeta \in \mathrm{Ker}[J]$. By construction, $\tilde{x}^k \in \mathrm{Im}[J] = \mathrm{Ker}[J^\top]^\perp = \mathrm{Ker}[J]^\perp$, hence $\zeta = 0$. Moreover, as $J$ is symmetrical and $\hat{x}(y) \in \mathrm{Im}[J]$, the quantity $J^+Jy$ coincides with $JJ^+y = \mathcal{R}^{\text{inv}}_{\hat{x}}(y)$ (by virtue of Remark 4), *i.e.,* the invariant re-fitting. Reminding the artifacts illustrated in Figure 3.(e), this is not satisfying.

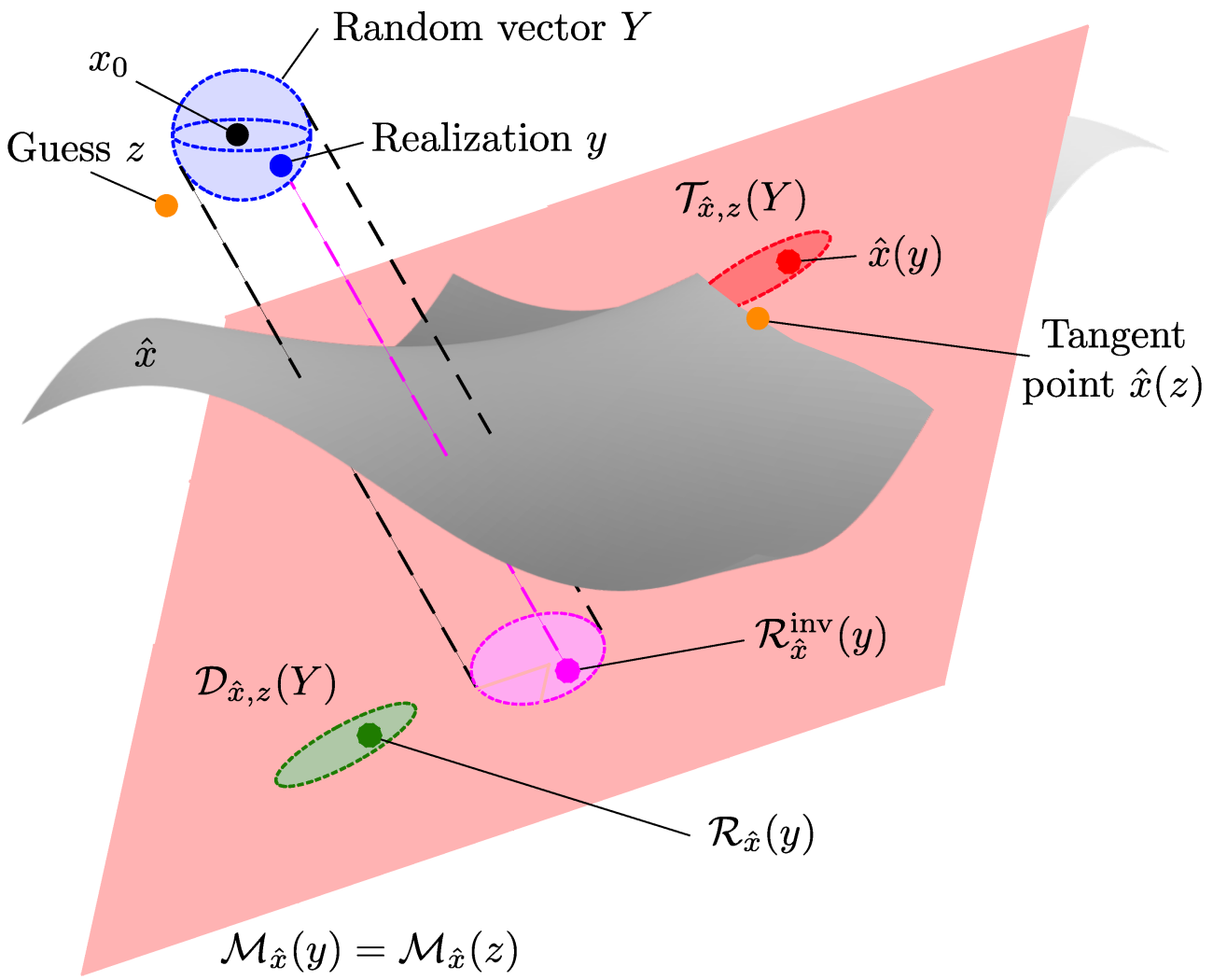

Fig. 5. *Geometrical illustration of the covariant re-fitting in a denoising problem of dimension $p = 3$. We assume that $\mathcal{M}_{\hat{x}}(z) = \mathcal{M}_{\hat{x}}(y)$ for the sake of clarity. The gray surface is the manifold modeling the evolution of $\hat{x}$ in an extended neighborhood of $y$. The light red affine plane is the model subspace tangent at $z$. The ellipses represent the positive-definite symmetric covariance matrices of some random vectors, as defined in Proposition 8 and Corollary 9.*

An interesting property of $\mathcal{R}_{\hat{x}}$ is the fact that it belongs to the model subspace of $\hat{x}$ as stated in the following proposition.

Proposition 17. *Let $y \mapsto \hat{x}(y)$ be an* a.e. *differentiable estimator. Then for almost all $y \in \mathbb{R}^n$, one has $\mathcal{R}_{\hat{x}}(y) \in \mathcal{M}_{\hat{x}}(y)$.*

*Proof.* As $\mathcal{M}_{\hat{x}}(y) = \hat{x}(y) + \mathrm{Im}\, J$ and $\rho J(y - \Phi\hat{x}(y)) \in \mathrm{Im}\, J$, the proposition holds. ∎

The case where $\Phi J$ is an orthogonal projector, leads to interesting properties that will be of interest regarding some estimators considered in Subsection 2.2.

Proposition 18. *Suppose that $\Phi J$ is an orthogonal projector. Then, $\mathcal{R}_{\hat{x}}(y) = \hat{x}(y) + J(y - \Phi\hat{x}(y))$, and, $\Phi\mathcal{R}_{\hat{x}}(y) = \Phi\mathcal{R}^{inv}_{\hat{x}}(y)$.*

*Proof.* By virtue of Proposition 7, $\rho = 1$ and then $\mathcal{R}_{\hat{x}}(y) = \hat{x}(y) + J(y - \Phi\hat{x}(y))$. The fact that $\Phi\mathcal{R}_{\hat{x}}(y) = \Phi\mathcal{R}^{\text{inv}}_{\hat{x}}(y)$ comes from the fact that $\Phi J(\Phi J)^+ = \Phi J$. ∎

The next proposition provides, when $J\Phi$ satisfies a fixed point formulation, an expression of $\mathcal{R}_{\hat{x}}(y)$ that will be useful for efficient computations of the re-fitting as discussed in Section 5, a notable example being the iso-TV regularization.

Proposition 19. *Assume that $J\Phi\hat{x}(y) = \hat{x}(y)$. Then, the covariant re-fitting reads $\mathcal{R}_{\hat{x}}(y) = (1 - \rho)\hat{x}(y) + \rho Jy$.*

*Proof.* We have $\mathcal{R}_{\hat{x}}(y) = \hat{x}(y) + \rho J(y - \Phi\hat{x}(y)) = \hat{x}(y) + \rho J y - \rho J \Phi \hat{x}(y)$, and since $J\Phi\hat{x}(y) = \hat{x}(y)$ by assumption, this concludes the proof. □

Interestingly, the next theorem shows that the condition $J\Phi\hat{x}(y) = \hat{x}(y)$ is met provided $\hat{x}(y)$ is solution of a variational problem with a 1-homogeneous regularizer.

THEOREM 20. *Let $\hat{x}(y)$ be the unique* a.e. *differentiable solution of*

$$\hat{x}(y) = \operatorname*{argmin}_x F(y - \Phi x) + G(x) \quad , \tag{35}$$

*with $F$, $G$ being convex and $G$ being 1-homogeneous. Then, $J\Phi\hat{x}(y) = \hat{x}(y)$* a.e. .

The proof of Theorem 20 is postponed to Appendix B.

The affine constrained least-squares, the $\ell_1$ synthesis, the $\ell_1 - \ell_2$ analysis, aniso-TV and iso-TV, are solutions of a variational problem with $F$ being differentiable and $G$ being 1-homogeneous. As a consequence, Theorem 20 shows that the aforementioned methods satisfy $J\Phi\hat{x}(y) = \hat{x}(y)$, and hence $\mathcal{R}_{\hat{x}}(y) = (1-\rho)\hat{x}(y) + \rho J y$.

**4.3. Examples of re-fitting procedures.** We now exemplify the previous definitions for the wide class of variational estimators introduced in Subsection 2.2.

**The affine constrained least-squares** have for Jacobian matrix $J = A(\Phi A)^+$. In this case, $\Phi J = \Phi A(\Phi A)^+$ is an orthogonal projector, $\rho = 1$ and the covariant re-fitting coincides with the invariant one and reads $\mathcal{R}_{\hat{x}}(y) = \mathcal{R}^{\text{inv}}_{\hat{x}}(y) = \hat{x}(y)$.

**The Tikhonov regularization** has for Jacobian $J = (\Phi^\top \Phi + \lambda \Gamma^\top \Gamma)^{-1} \Phi^\top$ and in this case $\rho$ depends on the residual $\delta$ and the re-fitting reads as the weighted sum $\mathcal{R}_{\hat{x}}(y) = (1+\rho)\hat{x}(y) - \rho J \Phi \hat{x}(y)$.

**The soft- and hard-thresholding**, used with $\Phi = \mathrm{Id}$, have the Jacobian matrix $J = \mathrm{Id}_{\mathcal{I}}$. As a consequence $\Phi J = \mathrm{Id}_{\mathcal{I}}$ is an orthogonal projection and the covariant re-fitting coincides with the invariant one, namely the hard-thresholding itself.

**The $\ell_1$ synthesis** has for Jacobian matrix $J = \mathrm{Id}_{\mathcal{I}}(\Phi_{\mathcal{I}})^+$ where $\Phi_{\mathcal{I}}$ has full column rank. As for the thresholding, $\Phi J = \Phi U(\Phi U)^+$ is an orthogonal projection and the covariant re-fitting reads $\mathcal{R}_{\hat{x}}(y) = \mathcal{R}^{\text{inv}}_{\hat{x}}(y)$.

**The $\ell_1$ analysis** has for Jacobian matrix $J = U(\Phi U)^+$. Again, $\Phi J = \Phi U(\Phi U)^+$ is an orthogonal projection and the covariant re-fitting reads $\mathcal{R}_{\hat{x}}(y) = \mathcal{R}^{\text{inv}}_{\hat{x}}(y)$.

**The $\ell_1 - \ell_2$ analysis** has the Jacobian operator given in Eq. (12), which applied to a vector $d \in \mathbb{R}^n$ is a solution of the following problem

$$Jd \in \operatorname*{argmin}_{x \,;\, \operatorname{supp}(\Gamma x) \subseteq \mathcal{I}} \tfrac{1}{2}\|\Phi x - d\|_2^2 + \tfrac{\lambda}{2}\omega(\Gamma x) \quad , \tag{36}$$

$$\text{where } \omega : z \in \mathbb{R}^{m \times b} \mapsto \sum_{i \in \mathcal{I}} \frac{1}{\|(\Gamma \hat{x}(y))_i\|_2} \left( \|z_i\|_2^2 - \left\langle z_i, \frac{(\Gamma \hat{x}(y))_i}{\|(\Gamma \hat{x}(y))_i\|_2} \right\rangle^2 \right) \quad .$$

Remark that $\omega(\Gamma J d) = 0$ only if $(\Gamma J d)_i$ is co-linear to $(\Gamma \hat{x}(y))_i$, for all $i \in \mathcal{I}$. For iso-TV, it means that the level lines of $Jd$ must be included in the ones of $\hat{x}(y)$. Moreover, by virtue of Theorem 20, one has $J\Phi\hat{x}(y) = \hat{x}(y)$ and hence $\mathcal{R}_{\hat{x}}(y) = (1-\rho)\hat{x}(y) + \rho J y$. As a consequence, unlike the invariant re-fitting of $\hat{x}(y)$, the covariant re-fitting is constrained to be faithful to the regularity of $\hat{x}(y)$, since it enforces the discontinuities of $Jd$ to be co-linear to $(\Gamma \hat{x}(y))_{\mathcal{I}}$. This is especially important where the iso-TV solution presents transitions with high curvature. Such an appealing behavior of the covariant re-fitting explains the results observed in Figure 2.(e) and Figure 3.(f).

**The non-local means** has an intricate Jacobian matrix. Nevertheless, its directional derivative has a simpler expression given, for any direction $d \in \mathbb{R}^n$, by

$$Jd = \frac{\sum_j \bar{w}'_{i,j} y_j + \sum \bar{w}_{i,j} d_j - \hat{x}(y)_i \sum_j \bar{w}'_{i,j}}{\sum_j \bar{w}_{i,j}} \quad \text{with} \quad \bar{w}'_{i,j} = \sum_k w'_{i-k,j-k} \ , \tag{37}$$

where $\hat{x}(y)$ is defined in Eq. (7) and

$$w'_{i,j} = 2\langle \mathcal{P}_i y - \mathcal{P}_j y, \, \mathcal{P}_i d - \mathcal{P}_j d \rangle \varphi' \left( \|\mathcal{P}_i y - \mathcal{P}_j y\|^2 \right) \ , \tag{38}$$

with $\varphi'$ the *a.e.* derivative of the kernel function $\varphi$. Subsequently, the covariant re-fitting is obtained from its general form with two steps, by computing first $\hat{x}(y)$, and applying next the Jacobian to the direction $d = y - \hat{x}(y)$.

**5. Covariant re-fitting in practice.** This section details the computation of CLEAR for standard algorithms. We first recall some properties of two different differentiation techniques that allow computing some image of $J$ jointly with $\hat{x}(y)$.

**5.1. Algorithmic differentiation.** Following [15], we consider restoration algorithms whose solutions $\hat{x}(y) = x^k$ are obtained via an iterative scheme of the form

$$\begin{cases} x^k & = \gamma(a^k) \ , \\ a^{k+1} & = \psi(a^k, y) \ . \end{cases} \tag{39}$$

Here, $a^k \in \mathcal{A}$ is a sequence of auxiliary variables, $\psi : \mathcal{A} \times \mathbb{R}^n \to \mathcal{A}$ is a fixed point operator in the sense that $a^k$ converges to $a^\star$, and $\gamma : \mathcal{A} \to \mathbb{R}^p$ is non-expansive (*i.e.*, $\|\gamma(a_1) - \gamma(a_2)\| \leqslant \|a_1 - a_2\|, \forall a_1, a_2 \in \mathcal{A}$) entailing $x^k$ converges to $x^\star = \gamma(a^\star)$.

As a result, for almost all $y$ and for any direction $d \in \mathbb{R}^n$, the directional derivatives $\mathcal{D}_x^k = J_{\hat{x}^k}(y) d$ and $\mathcal{D}_a^k = J_{a^k}(y) d$ can be jointly obtained with $x^k$ and $a^k$ as

$$\begin{cases} x^k & = \gamma(a^k) \ , \\ a^{k+1} & = \psi(a^k, y) \ , \\ \mathcal{D}_x^k & = \Gamma_a \mathcal{D}_a^k \ , \\ \mathcal{D}_a^{k+1} & = \Psi_a \mathcal{D}_a^k + \Psi_y d \ , \end{cases} \tag{40}$$

where $\Gamma_a = \left.\frac{\partial \gamma(a)}{\partial a}\right|_{a^k}$, $\Psi_a = \left.\frac{\partial \psi(a,y)}{\partial a}\right|_{a^k}$ and $\Psi_y = \left.\frac{\partial \psi(a^k,y)}{\partial y}\right|_{y}$. Interestingly, in all considered cases, the cost of evaluating $\Gamma_a$, $\Psi_a$ and $\Psi_y$ is about the same as the one of evaluating $\gamma$ and $\psi$. As a result, the complexity of (40) is of about twice the complexity of (39). In practice, $\Gamma_a$, $\Psi_a$ and $\Psi_y$ can be implemented either thanks to their closed form expression or in a black box manner using automatic differentiation. The later has been well studied and we refer to [24, 31] for a comprehensive study.

**5.2. Finite difference based differentiation.** Another strategy is to approximate directional derivatives by finite differences, for any $d \in \mathbb{R}^n$ and $\varepsilon > 0$, as

$$J_{\hat{x}}(y) d \approx \frac{\hat{x}(y + \varepsilon d) - \hat{x}(y)}{\varepsilon} \ . \tag{41}$$

As a result, the complexity of evaluating (41) is also twice the complexity of (39) since $\hat{x}$ must be evaluated at both $y$ and $y + \varepsilon d$. The main advantage of this method is that $\hat{x}$ can be used as a black box, *i.e.*, without any knowledge on the underlying algorithm that provides $\hat{x}(y)$. For $\varepsilon$ small enough, it performs as well as the approach described in (40) (with $\hat{x}(y) = x^k$) that requires the knowledge of the derivatives. Indeed, if $y \mapsto \hat{x}(y)$ is Lipschitz-continuous, then (41) converges to (40) when $\varepsilon \to 0$ (by virtue of Rademacher's theorem and [21, Theorem 1-2, Section 6.2]). This implies that the value $\varepsilon$ can be chosen as small as possible (up to machine precision) yielding an accurate approximation of $J_{\hat{x}}(y)d$. This finite difference approach has been used in many fields, and notably for risk estimation, see *e.g.,* [53, 40, 35]. We will apply this black box strategy on state-of-the-art algorithms in Subsection 6.5.

**5.3. Two-step computation for the general case.** In the most general case, the computation of the covariant re-fitting, given by

$$(42)\quad \mathcal{R}_{\hat{x}}(y) = \hat{x}(y) + \rho J(y - \Phi\hat{x}(y)) \quad \text{with} \quad \rho = \frac{\langle \Phi J\delta,\, \delta\rangle}{\|\Phi J\delta\|_2^2} \quad \text{and} \quad \delta = y - \Phi\hat{x}(y)\ ,$$

requires to evaluate sequentially $\hat{x}(y)$ and $J(y - \Phi\hat{x}(y))$.

With finite difference differentiation, two steps are required. First $\hat{x}(y)$ must be computed with the original algorithm and next $J(y - \Phi\hat{x}(y))$ is obtained by finite difference (41) on the direction of the residual $d = y - \Phi\hat{x}(y)$. Once $J(y - \Phi\hat{x}(y))$ is computed, $\rho$ can be evaluated and subsequently (42). The overall complexity is about twice that of the original algorithm producing $\hat{x}(y)$.

With algorithmic differentiation, as $J(y - \Phi\hat{x}(y))$ depends on $\hat{x}(y)$, the original iterative scheme (39) must be run first. In the second step, $J(y - \Phi\hat{x}(y))$ is obtained with the differentiated version (40) on the direction of the residual $d = y - \Phi\hat{x}(y)$. As a result, $\hat{x}(y)$ is computed twice, first by (39), next by (40). It leads to an overall complexity about three times the one of the original algorithm. Nevertheless, in several cases, one can avoid the first step by running (40) only once.

**5.4. One-step computation for specific cases.** When $\hat{x}(y)$ fulfills the assumption $J\Phi\hat{x}(y) = \hat{x}(y)$ of Proposition 19, the covariant re-fitting reads as

$$(43)\qquad \mathcal{R}_{\hat{x}}(y) = (1-\rho)\hat{x}(y) + \rho J y \quad \text{with} \quad \rho = \frac{\langle \Phi(Jy - \hat{x}(y)),\, y - \Phi\hat{x}(y)\rangle}{\|\Phi(Jy - \hat{x}(y))\|_2^2}\ .$$

The computations of $\hat{x}(y)$ and $Jy$ are then sufficient to compute the re-fitting $\mathcal{R}_{\hat{x}}(y)$. As a result, in the case of algorithmic differentiation, (40) can be run once to get $\mathcal{R}_{\hat{x}}(y)$ since using $d = y$ provides directly $\hat{x}(y)$, $Jy$ and subsequently $\rho$. Compared to the two step approach, the complexity of the re-fitting reduces to about twice the one of the original step from (39). Recall, that the condition $J\Phi\hat{x}(y) = \hat{x}(y)$ is met for the Lasso, the Generalize Lasso, aniso-TV and iso-TV, hence they can be re-enhanced with a complexity being twice the one of their original algorithm.

**5.5. Example on a primal-dual solver for $\ell_1$ analysis.** In this section we instantiate Algorithm (40) to the case of the primal-dual sequence of [7]. By dualizing the $\ell_1$ analysis norm $x \mapsto \lambda\|\Gamma x\|_1$, the primal problem (4) can be reformulated, with $x^\star = \hat{x}(y)$, as the following saddle-point problem

$$(44)\qquad (z^\star, x^\star) \in \arg\max_{z\in\mathbb{R}^m} \min_{x\in\mathbb{R}^p} \frac{1}{2}\|\Phi x - y\|_2^2 + \langle \Gamma x, z\rangle - \iota_{\mathcal{B}_\lambda}(z)\ ,$$

where $z^\star \in \mathbb{R}^m$ is the dual variable, and $\mathcal{B}_\lambda = \{z \in \mathbb{R}^m \ :\ \|z\|_\infty \leqslant \lambda\}$ is the $\ell_\infty$ ball.

Problem (44) can be efficiently solved using the primal-dual algorithm of [7]. By taking $\sigma\tau < 1/\|\Gamma\|_2^2$, $\theta \in [0,1]$ and initializing (for instance,) $x^0 = v^0 = 0 \in \mathbb{R}^p$,

$z^0 = 0 \in \mathbb{R}^m$, the algorithm reads

$$(45)\qquad \begin{cases} z^{k+1} &= \Pi_{\mathcal{B}_\lambda}(z^k + \sigma \Gamma v^k) \ , \\ x^{k+1} &= (\mathrm{Id} + \tau \Phi^\top \Phi)^{-1} \left( x^k + \tau(\Phi^\top y - \Gamma^\top z^{k+1}) \right) \ , \\ v^{k+1} &= x^{k+1} + \theta(x^{k+1} - x^k) \ , \end{cases}$$

where the projection of $z$ over $\mathcal{B}_\lambda$ is done component-wise as

$$(46)\qquad \Pi_{\mathcal{B}_\lambda}(z)_i = \begin{cases} z_i & \text{if} \quad |z_i| \leqslant \lambda \ , \\ \lambda \ \mathrm{sign}(z_i) & \text{otherwise} \ . \end{cases}$$

The sequence $x^k$ converges to a solution $x^\star$ of the $\ell_1$ analysis problem [7].

It is easy to check that the primal-dual sequence defined in (45) can be written in the general form considered in (39), see for instance [15]. As a result, we can use the algorithmic differentiation based strategy described by (40) as follows: for the initialization $\tilde{x}^0 = \tilde{v}^0 = 0 \in \mathbb{R}^p$, $\tilde{z}^0 = 0 \in \mathbb{R}^m$, and for $\beta = 0$, as

$$(47)\qquad \begin{cases} z^{k+1} &= \Pi_{\mathcal{B}_\lambda}(z^k + \sigma \Gamma v^k) \ , \\ x^{k+1} &= (\mathrm{Id} + \tau \Phi^\top \Phi)^{-1} \left( x^k + \tau(\Phi^\top y - \Gamma^\top z^{k+1}) \right) \ , \\ v^{k+1} &= x^{k+1} + \theta(x^{k+1} - x^k) \ , \\ \tilde{z}^{k+1} &= \Pi_{z^k + \sigma \Gamma v^k}(\tilde{z}^k + \sigma \Gamma \tilde{v}^k) \ , \\ \tilde{x}^{k+1} &= (\mathrm{Id} + \tau \Phi^\top \Phi)^{-1} \left( \tilde{x}^k + \tau(\Phi^\top y - \Gamma^\top \tilde{z}^{k+1}) \right) \ , \\ \tilde{v}^{k+1} &= \tilde{x}^{k+1} + \theta(\tilde{x}^{k+1} - \tilde{x}^k) \ , \end{cases}$$

$$\text{where} \quad \Pi_z(\tilde{z})_i = \begin{cases} \tilde{z}_i & \text{if} \quad |z_i| \leqslant \lambda + \beta \ , \\ 0 & \text{otherwise} \ . \end{cases}$$

Recall that the re-fitting is $\mathcal{R}_{x^k}(y) = \tilde{x}^k$, since $J\Phi$ is an orthogonal projector.

Remark that the algorithmic differentiation of (45) is exactly (47) for $\beta = 0$, hence, $\tilde{x}^k = \mathcal{R}_{x^k}(y)$. However, if one wants to guarantee the convergence of the sequence $\tilde{x}^k$ towards $\mathcal{R}_{\hat{x}}(y)$, one needs a small $\beta > 0$ as shown in the next theorem. In practice, $\beta$ can be chosen as the smallest available positive floating number.

Theorem 21. *Assume that $x^\star$ satisfies (5) with $\Phi U$ full-column rank[2]. Let $\alpha > 0$ be the minimum non zero value[3] of $|\Gamma x^\star|_i$ for all $i \in [m]$. Choose $\beta$ such that $\alpha\sigma > \beta > 0$. Then, the sequence $\tilde{x}^k = \mathcal{R}_{x^k}(y)$ defined in (47) converges to the re-fitting $\mathcal{R}_{\hat{x}}(y)$ of $\hat{x}(y) = x^\star$.*

The proof of this theorem is postponed to Appendix C.

A similar result was obtained in [13] when solving the $\ell_1$ analysis problem (4) with the Douglas-Rachford splitting algorithm described in [17, 8].

**5.6. Example for a primal-dual solver for $\ell_1 - \ell_2$ analysis.** The algorithm for the $\ell_1 - \ell_2$ analysis regularization can be derived with the exact same considerations as for the $\ell_1$ analysis case. The only difference in the application of the primal dual algorithm comes from the non linear operation (46) that now reads, for $z \in \mathbb{R}^{m \times b}$, as

$$(48)\qquad \Pi^{\mathrm{iso}}_{\mathcal{B}_\lambda}(z)_i = \begin{cases} z_i & \text{if} \quad \|z_i\|_2 \leqslant \lambda \ , \\ \lambda \frac{z_i}{\|z_i\|_2} & \text{otherwise} \ . \end{cases}$$

[2]This could be enforced as shown in [48].

[3]If $|\Gamma x^\star|_i = 0$ for all $i \in [m]$, the result remains true for any $\alpha > 0$.

**Algorithm** Non-local means [3] and its directional derivative

**Inputs:** noisy image $y$, direction $d$, noise standard-deviation $\sigma$
**Parameters:** half search window width $s$, half patch width $b$, kernel function $\varphi$
**Outputs:** $x^\star = \hat{x}(y)$ and $\tilde{x} = J_{\hat{x}}(y)d$

```
Initialize  W ← φ(2σ²(2b+1)²) 1_{p1×p2}          (add weights for the central pixels [39])
Initialize  W_y ← φ(2σ²(2b+1)²) y                 (accumulators for the weighted sum)
Initialize  W' ← 0_{p1×p2}
Initialize  W'_y ← φ(2σ²(2b+1)²) d
for k ∈ [−s, s]² \ {0, 0} do
    Compute  e ← [ (y − S_k(y))² ] ⋆ κ           (error between each k shifted patches [11, 12])
    Compute  w ← φ(e) ⋆ κ                         (contribution for each patch of its k shift)
    Update   W ← W + w                            (add weights at each position)
    Update   W_y ← W_y + w S_k(y)                 (add contribution of each k shifted patches)

    Compute  e' ← [ 2(y − S_k(y))(d − S_k(d)) ] ⋆ κ
    Compute  w' ← [ e' φ'(e) ] ⋆ κ
    Update   W' ← W' + w'
    Update   W'_y ← W'_y + w' S_k(f) + w S_k(d)
end for
Compute  x⋆ ← W_y / W                             (weighted mean)
Compute  x̃ ← (W'_y − W' x⋆) / W
```

FIG. 6. *Pseudo-algorithm for the computation of the non-local means and its Jacobian in a direction $d$. All arithmetic operations are element wise, $S_k$ is the operator that shift all pixels in the direction $k$, $\star$ is the discrete convolution operator, and $\kappa \in \mathbb{R}^{p_1 \times p_2}$ is such that $\kappa_{i,j} = 1$ if $(i, j) \in [-b, b]^2$, 0 otherwise.*

It follows that the algorithmic differentiation strategy reads as

$$\text{(49)} \quad \left\{ \begin{array}{ll} z^{k+1} & = \Pi^{\text{iso}}_{\mathcal{B}_\lambda}(z^k + \sigma\Gamma v^k) \ , \\ x^{k+1} & = (\mathrm{Id} + \tau\Phi^\top\Phi)^{-1}\left(x^k + \tau(\Phi^\top y - \Gamma^\top(z^{k+1}))\right) \ , \\ v^{k+1} & = x^{k+1} + \theta(x^{k+1} - x^k) \ , \\ \tilde{z}^{k+1} & = \Pi^{\text{iso}}_{z^k + \sigma\Gamma v^k}(\tilde{z}^k + \sigma\Gamma\tilde{v}^k) \ , \\ \tilde{x}^{k+1} & = (\mathrm{Id} + \tau\Phi^\top\Phi)^{-1}\left(\tilde{x}^k + \tau(\Phi^\top y - \Gamma^\top \tilde{z}^{k+1})\right) , \\ \tilde{v}^{k+1} & = \tilde{x}^{k+1} + \theta(\tilde{x}^{k+1} - \tilde{x}^k) \ , \end{array} \right.$$

$$\text{where} \quad \Pi^{\text{iso}}_z(\tilde{z})_i = \begin{cases} \tilde{z}_i & \text{if} \quad \|z_i\|_2 \leqslant \lambda + \beta \ , \\ \frac{\lambda}{\|z_i\|_2}\left(\tilde{z}_i - \left\langle \tilde{z}_i, \frac{z_i}{\|z_i\|_2} \right\rangle \frac{z_i}{\|z_i\|_2}\right) & \text{otherwise} \ . \end{cases}$$

Unlike the $\ell_1$ case, the re-fitted solution is not $\tilde{x}^k$ itself, but following Subsection 5.4, it can be obtained at iteration $k$ as

$$\text{(50)} \qquad \mathcal{R}_{x^k}(y) = (1 - \rho)x^k + \rho\tilde{x}^k \quad \text{with} \quad \rho = \frac{\langle \Phi(\tilde{x}^k - x^k),\, y - \Phi x^k \rangle}{\|\Phi(\tilde{x}^k - x^k)\|_2^2} \ .$$

**5.7. Example for the non-local means.** In this section, we specify the update rule (40) to an acceleration of the non-local means inspired from [11, 12]. We use the procedure of [39] to correctly handle the central pixel. Again, one can check that this implementation can be written in the general form considered in the update rule (39), where the fixed point solution is obtained directly at the first iteration.

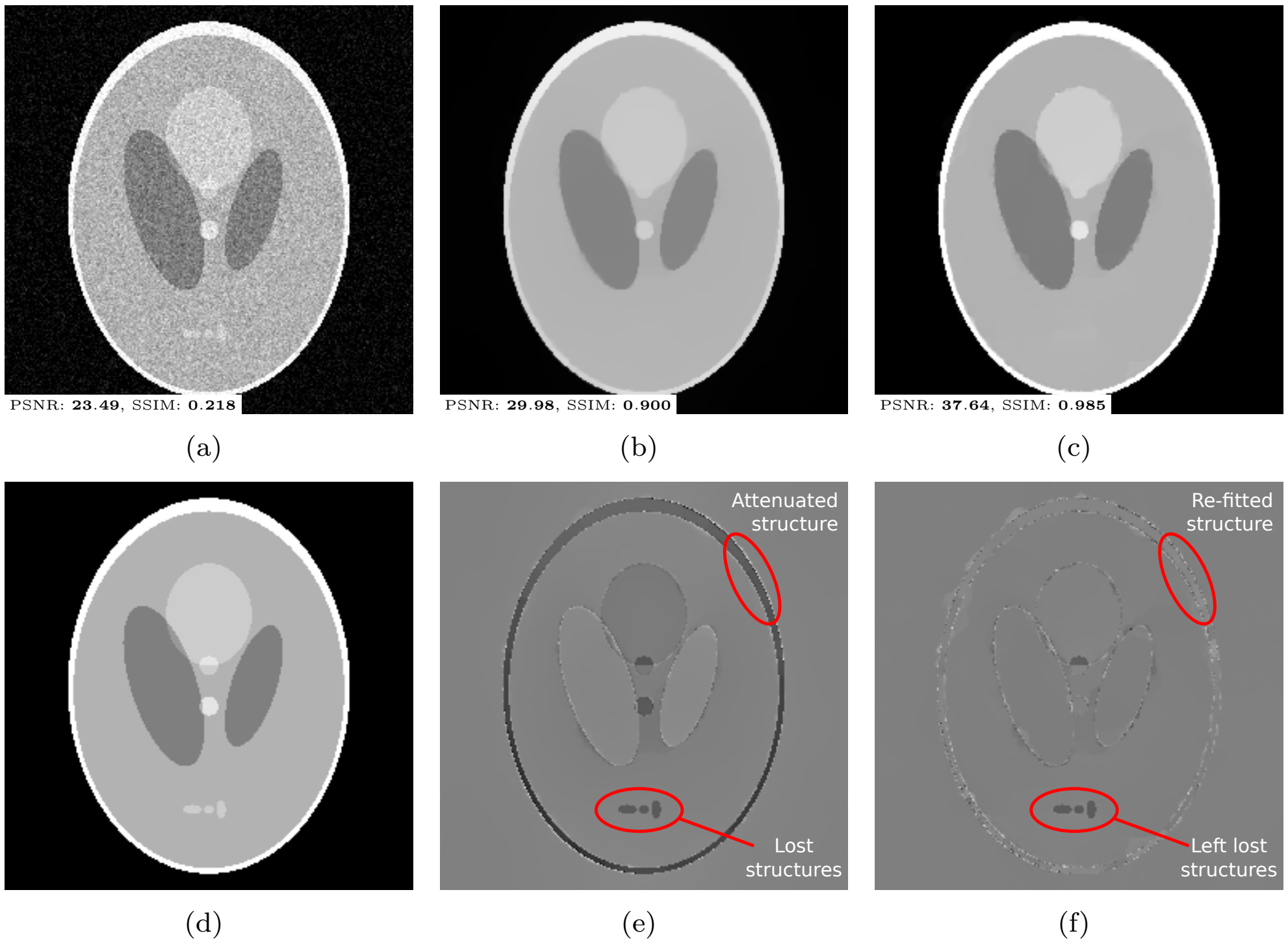

FIG. 7. *(a) Noisy* $y = x_0 + w$ *and (d) noise-free* $x_0$*. (b) Iso-TV* $\hat{x}(y)$ *with* $\lambda = 3$ *and (e) the residual* $\hat{x}(y) - x_0$*. (c) Our re-fitting* $\mathcal{R}_{\hat{x}}(y)$ *and (f) the residual* $\mathcal{R}_{\hat{x}}(y) - x_0$*.*

The pseudo-code obtained with the algorithmic differentiation scheme, as described in (40) is given in Figure 6. All variables with suffix $'$ correspond to the directional derivative obtained by using the chain rule on the original variables. This fast computation relies on that all convolutions can be computed with integral tables, leading to a global complexity in $O(s^2 n)$, for both the computation of the estimator $\hat{x}(y)$ and its directional derivative $Jd$. Recall that the covariant re-fitting is obtained from its general form with two steps, by computing first $\hat{x}(y)$, and applying next the proposed pseudo-code in the direction $d = y - \hat{x}(y)$.

## 6. Numerical experiments and comparisons with related approaches.

Here, we first give illustrations of our CLEAR method on toy image restoration problems. Then, through quantitative results in terms of PSNR[4] and SSIM [51], we numerically evaluate the re-fitting, discuss its benefit and limitations in several scenarios and compare this method with popular approaches from the literature.

### 6.1. Denoising with isotropic total-variation (iso-TV).

Figure 7 gives an illustration of our covariant re-fitting of the 2D iso-TV, where $\lambda$ has been chosen large enough in order to highlight the behavior of the re-fitting. We apply it for the denoising (*i.e.,* $\Phi = \mathrm{Id}$) of an $8bits$ piece-wise constant image damaged by AWGN with standard deviation $\sigma = 20$, known as the *Shepp-Logan phantom*. As discussed earlier, iso-TV introduces a significant loss of contrast [42], typically for thin detailed structures, which are re-enhanced in our result.

[4] $\mathrm{PSNR} = 10 \log_{10} 255^2 / \frac{1}{p} \|\hat{x}(y) - x_0\|_2^2$ for an image ranging on $[0, 255]$.

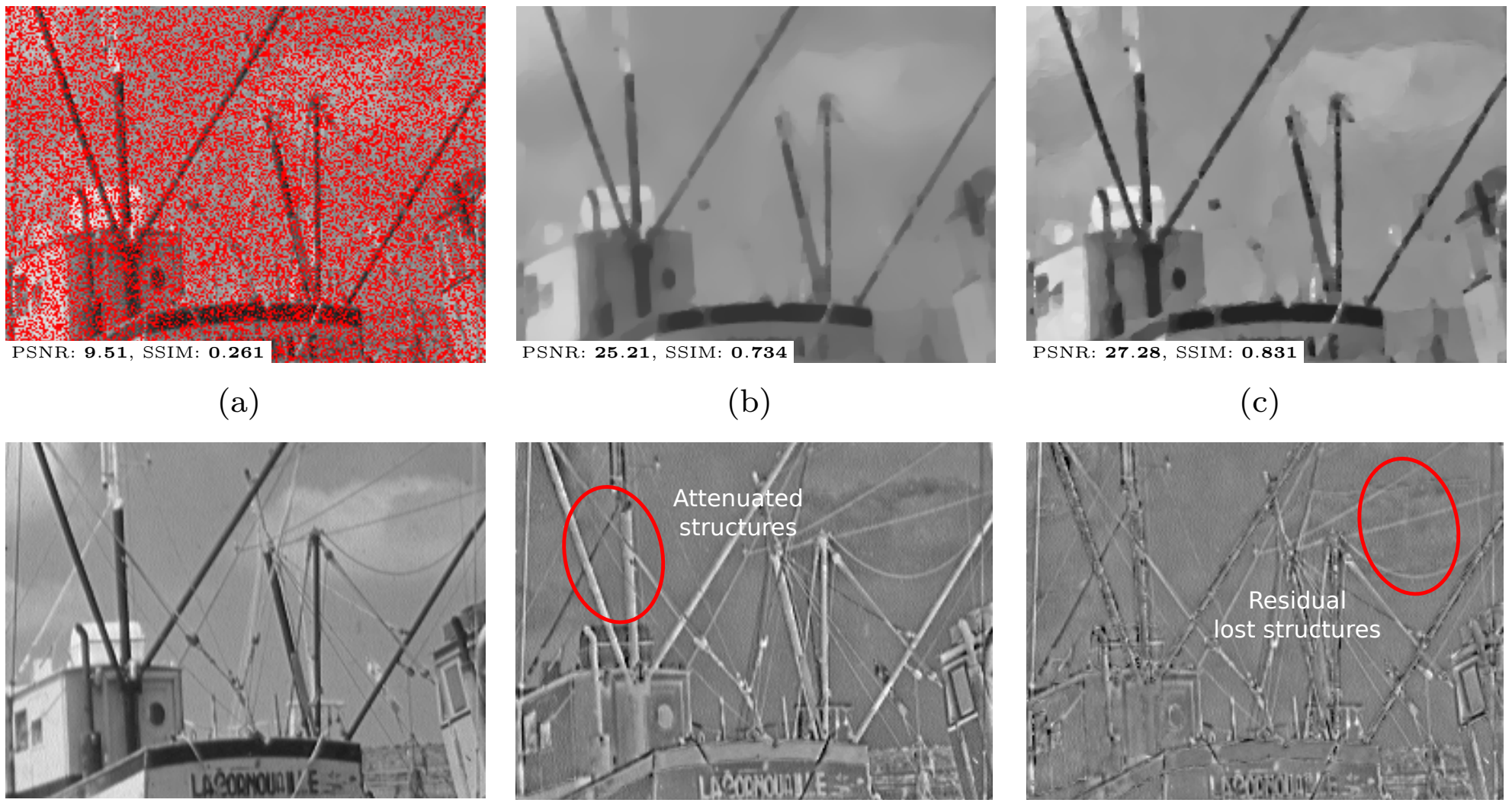

(a) (b) (c)

(d) (e) (f)

FIG. 8. *(a) Partial and noisy $y = \Phi x_0 + w$ (red indicates missing pixels) and (d) noise-free $x_0$. (b) Iso-TV $\hat{x}(y)$ with $\lambda = 3$ and (e) the residual $\hat{x}(y) - x_0$. (c) Our re-fitting $\mathcal{R}_{\hat{x}}(y)$ and (f) the residual $\mathcal{R}_{\hat{x}}(y) - x_0$.*

The residuals $\hat{x}(y) - x_0$ and $\mathcal{R}_{\hat{x}}(y) - x_0$ highlight that our re-fitting technique re-enhances efficiently the attenuated structure while leaving the lost structures unchanged. Nevertheless, after re-fitting, some small residuals around the edges appear. In fact, in the vicinity of edges, iso-TV finds (barely visible) discontinuities that are not in accordance with the underlying image. This creates an overload of small constant regions. When re-fitting is performed, such regions are re-fitted to the noisy data, and they become barely visible artifacts. In other words, the re-fitting has re-enforced the presence of a modeling problem, resulting to an increase of residual variance, that iso-TV had originally compensated by attenuating the amplitudes.

**6.2. De-masking with isotropic total-variation (iso-TV).** Figure 8 gives another illustration of our covariant re-fitting of the 2D iso-TV used for the restoration of an approximate 8*bits* piece-wise constant image damaged by AWGN with $\sigma = 20$, known as *Boat*. The operator $\Phi$ is a random mask removing 40% of pixels. Again, iso-TV introduces a significant loss of contrast, typically for thin details such as the contours of the objects, which are re-enhanced in our re-fitting result.

Inspecting the map of residuals in Figure 8.(e)-(f) illustrates that our re-fitting technique eliminates most of the bias to the price of a small variance increase. This is clearer when looking at the mast and the ropes. While the mast was preserved by iso-TV and re-enhanced by our re-fitting, the ropes remain lost for both methods.

**6.3. Denoising with non-local means.** Figure 9 gives an illustration of our re-fitting procedure for the non-local means algorithm used in a denoising problem of the 8*bits* image *Pirate*, enjoying many repetitive patterns and damaged by AWGN with $\sigma = 20$. We choose a regularizing kernel $\varphi(\cdot) = \exp(\cdot/h)$, $h > 0$ that leads to strong smoothing in order to highlight the behavior of the re-fitting. Our re-fitting technique provides favorable results: many details are enhanced compared to the standard method. This reveals that the non-local means is actually able to well capture the repetitions of many patterns but this information is not used properly to create a satisfying result. The re-fitting produces a sharper result, by enforcing the correct use of all the structures identified by patch comparisons. The maps of residuals Figure 9.(e)-(f) highlight that our re-fitting technique suppresses efficiently this dull effect while it preserves the model originally captured by patch redundancy. Again the suppression of this phenomenon is counter balanced by an increase of the residual variance, prominent where the local patch redundancy assumption is violated.

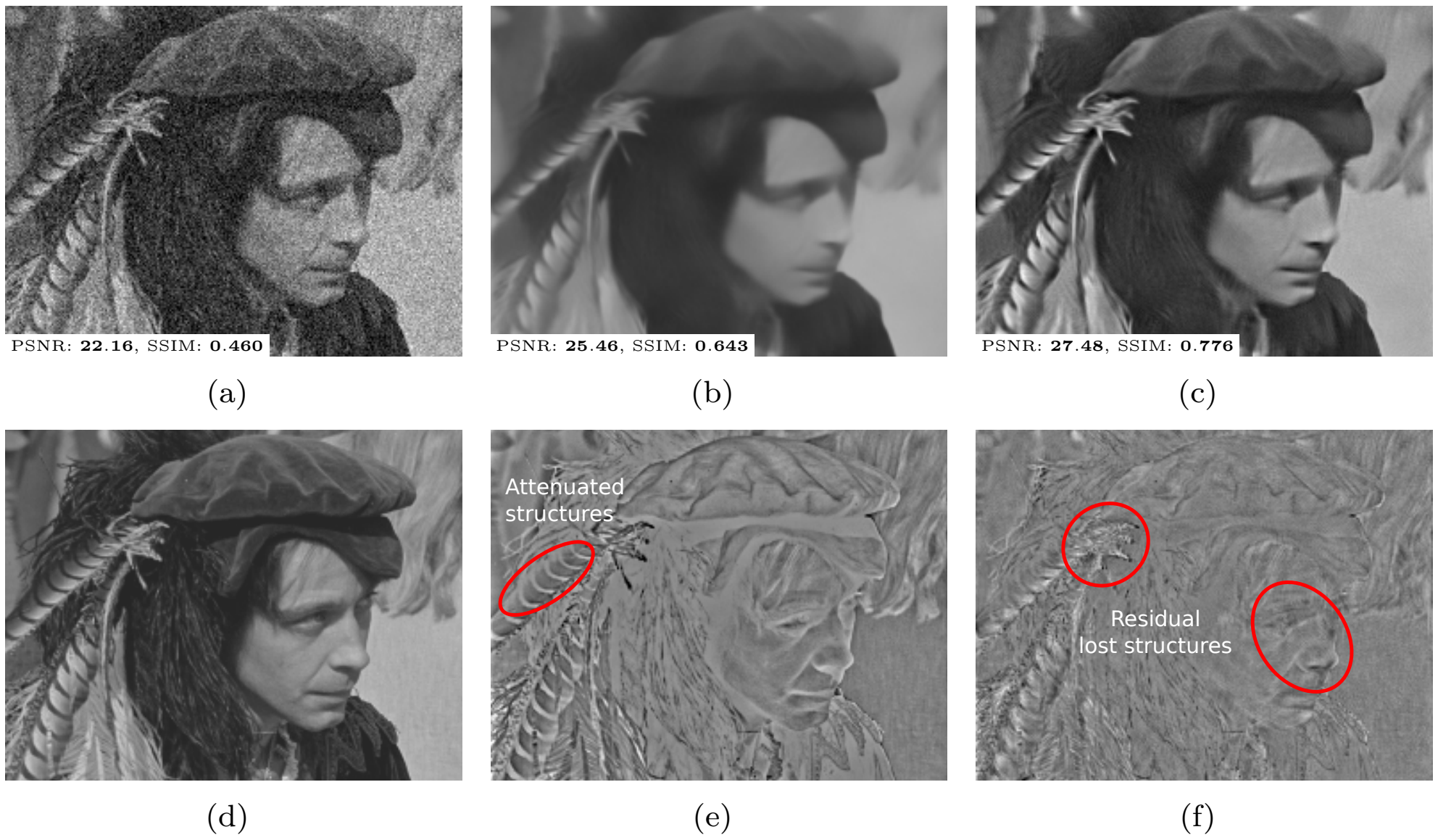

Fig. 9. *(a) Noisy* $y = x_0 + w$ *and (d) noise-free* $x_0$*. (b) Non-local means* $\hat{x}(y)$ *with* $h = 3$ *and (e) the residual* $\hat{x}(y) - x_0$*. (c) Our re-fitting* $\mathcal{R}_{\hat{x}}(y)$ *and (f) the residual* $\mathcal{R}_{\hat{x}}(y) - x_0$*.*

In these examples, the overall residual norm is clearly reduced by the re-fitting because the amount of reduced bias surpasses the increase of residual variance. This favorable behavior depends on the internal parameters of the original estimator acting on the bias-variance trade-off, as we investigate next.

**6.4. A bias-variance analysis for the covariant re-fitting.** Previous experiments have revealed that while CLEAR tends to reduce the bias, it increases (as expected) the residual variance. It is therefore important to understand under which conditions the bias-variance trade-off is in favor of our re-fitting technique.

Figure 10 illustrates the evolution of performance, measured in terms of mean squared error (MSE), of both aniso-TV and its re-fitting version as a function of the regularization parameter $\lambda$. Two images are considered: *Cameraman*, an approximate piece-wise constant image (top), and a truly piece-wise constant image (bottom).

This experiment highlights that optimal results for both approaches are not reached at the same $\lambda$ value. Visual inspection of the optima shows that due to the bias, the optimal solution of aniso-TV is reached for a $\lambda$ value promoting a model subspace that is not in accordance with the underlying signal: typically the presence of an overload of (barely visible) transitions in homogeneous areas. These transitions become clear when looking at the re-fitted version where each small region is re-fitted on the noisy data, revealing an excessive residual variance. Conversely, the optimal $\lambda$ value for the re-fitting seems to retrieve the correct model, *i.e.,* with transitions that are closely in accordance with the underlying signal. Comparing their relative performance, when both are used at their own optimal $\lambda$, reveals that our re-fitting brings a significant improvement if the underlying image is in fact piece-wise constant.

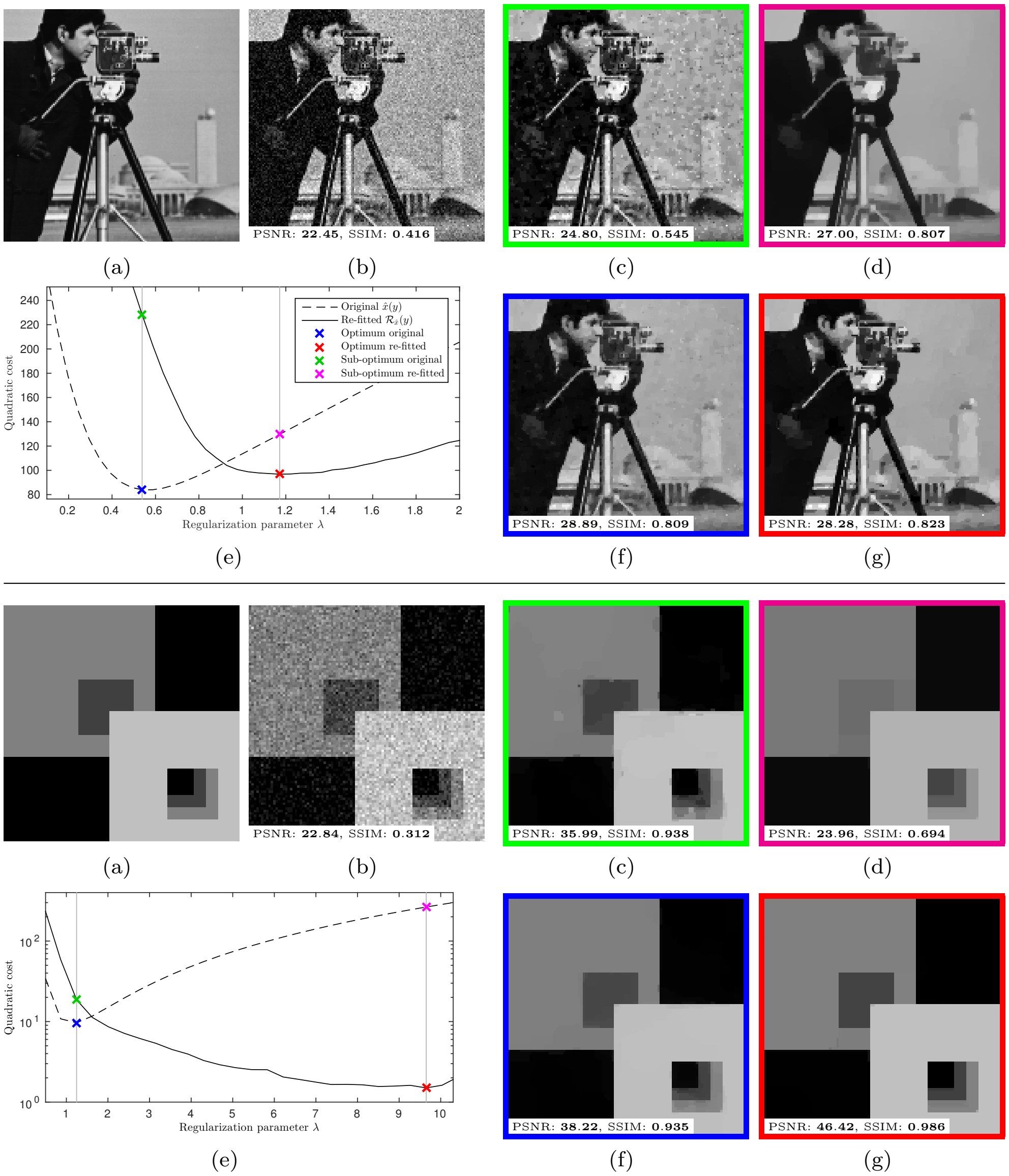

FIG. 10. *Experiment with aniso-TV: (top) poorly piece-wise constant case. (bottom) pure piece-wise constant case. (a) Noise-free $x_0$. (b) Noisy $y = x_0 + w$. (e) MSE of $\hat{x}(y)$ and its re-fitting* $\mathcal{R}_{\hat{x}}(y)$ w.r.t. $\lambda$. *Two values of $\lambda$ are selected corresponding to (c) re-fitting for a sub-optimal $\lambda$, (d) original for a sub-optimal $\lambda$, (f) original for the optimal $\lambda$, (g) re-fitting for the optimal $\lambda$.*

Figure 11 provides a similar illustration of the evolution of performance for the non-local means and its re-fitted version as a function of the smoothing parameter $h$ of the kernel function $\varphi(\cdot) = \exp(\cdot/h)$. Two images are considered: *Lady*, a crude approximation of an image with redundant patterns (top), and *Fingerprint*, an image with many redundant patterns (bottom). Similar conclusions can be made from this experiment. In particular, comparing their relative performance, when both are used at their own optimal $h$ value, seems to demonstrate that the re-fitting brings an improvement when most patches of the underlying image are redundant.

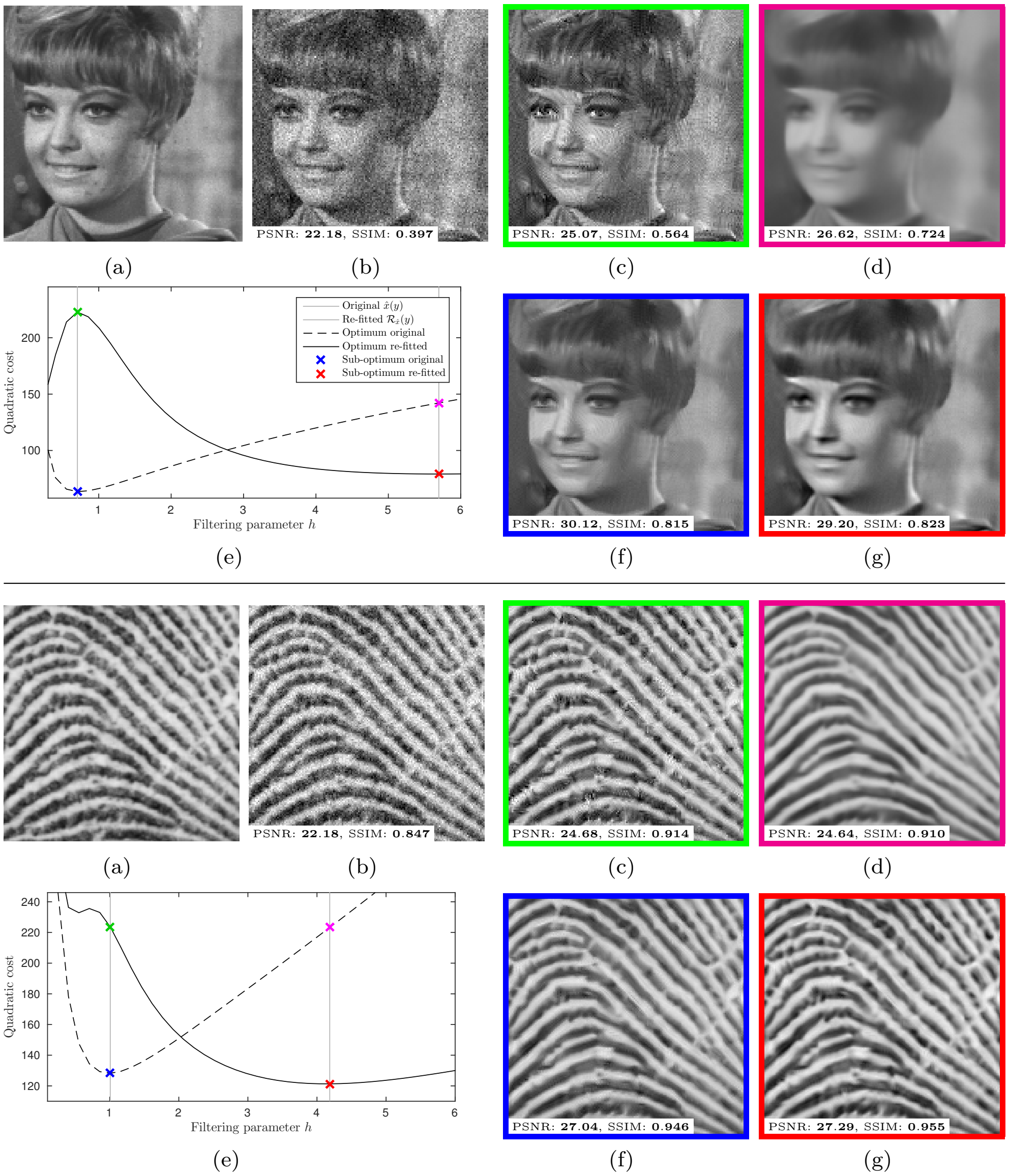

FIG. 11. *Experiment with non-local means: (top) moderate patch redundancy case. (bottom) high patch redundancy case. (a) Noise-free $x_0$. (b) Noisy $y = x_0 + w$. (e) MSE of $\hat{x}(y)$ and its re-fitting $\mathcal{R}_{\hat{x}}(y)$* w.r.t. *$h$. Two values of $h$ are selected corresponding to (c) re-fitting for a sub-optimal $h$, (d) original for a sub-optimal $h$, (f) original for the optimal $h$, (g) re-fitting for the optimal $h$.*

While it is difficult to make a general statement, we can reasonably claim from these experiments that re-fitting is all the more relevant in terms of MSE than the underlying image $x_0$ is in agreement with the retrieved subspace model $\mathcal{M}_{\hat{x}}(y)$. In other words, re-fitting is safe when the original restoration technique was chosen appropriately *w.r.t.* the underlying image of interest. Beyond MSE performance, the re-fitted results at their optimal parameters choice might nevertheless be preferable (as assessed by the SSIM values): intensities and contrasts being recovered better.

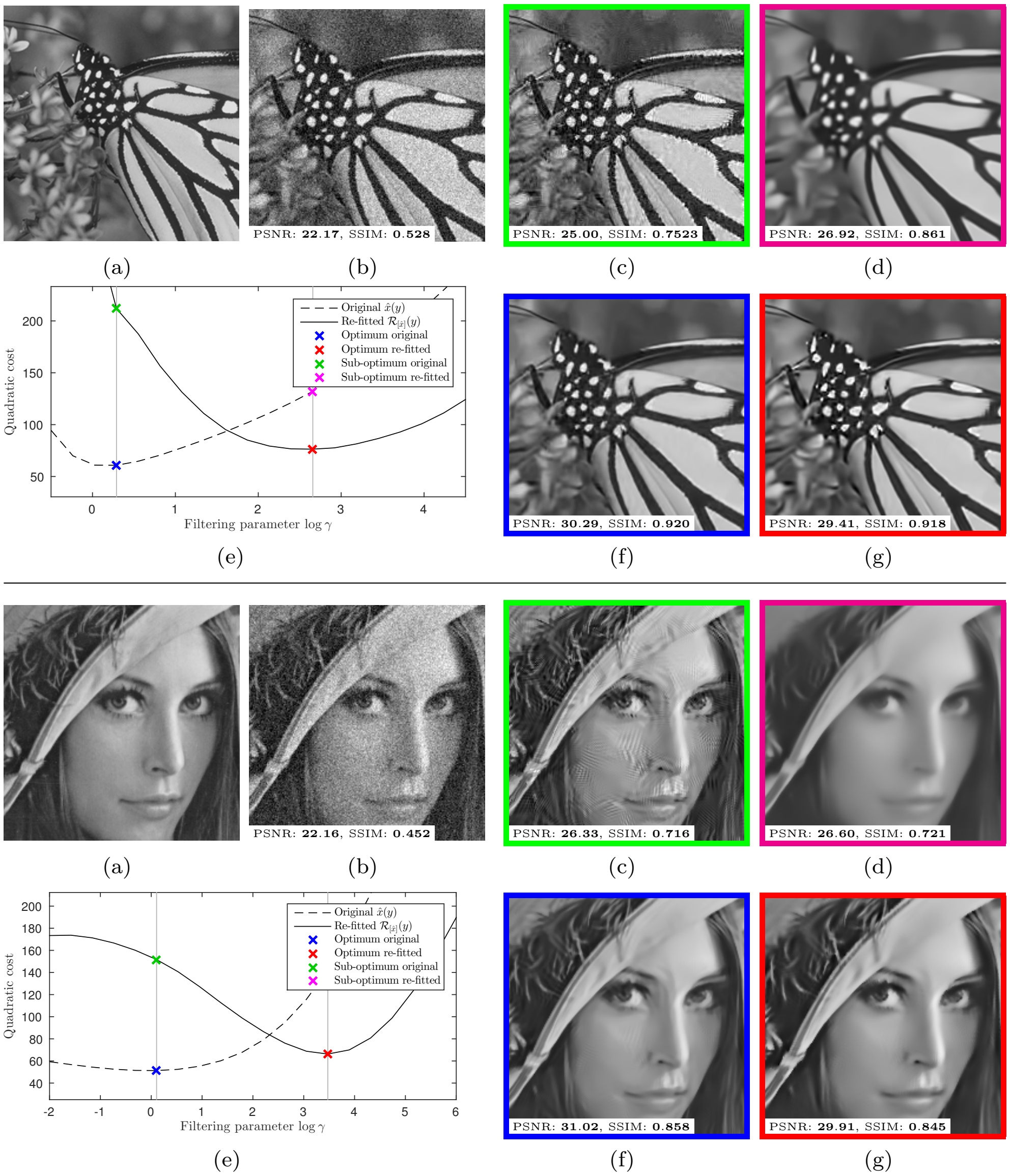

FIG. 12. *(top) Experiments with BM3D [10, 27] and with DDID [26] (bottom). (a) Noise-free $x_0$. (b) Noisy $y = x_0 + w$. (e) MSE of $\hat{x}(y)$ and its re-fitting $\mathcal{R}_{\hat{x}}(y)$* w.r.t. *$\gamma$. Two values of $\gamma$ are selected corresponding to (c) re-fitting for a sub-optimal $\gamma$, (d) original for a sub-optimal $\gamma$, (f) original for its optimal $\gamma$, (g) re-fitting for its optimal $\gamma$.*

**6.5. Behavior on more sophisticated filters.** We focus here on two filters BM3D [10] and DDID [26] reaching state-of-the-art results in denoising. As mentioned in the previous paragraph, re-fitting is all the more relevant that the first estimate is obtained with high smoothing strength. It is thus important to compare the original filter with its re-fitted version for varying smoothing parameters. For simplicity, we have considered only one smoothing parameter $\gamma > 0$ for these two algorithms.

For DDID, we use the authors implementation, and choose to multiply by $\gamma$ the two inner parameters $\gamma_f$ and $\gamma_r$ [26]. For BM3D, we use the implementation of [27], and choose to multiply by $\gamma$ the inner parameters $\lambda_{3D}^{\mathbf{hard}}$ (in the hard thresholding step) and $\sigma$ (in the Wiener filtering step). Unlike previous experiments, the re-fitted results are obtained here by finite difference as discussed in Subsection 5.2.

Figure 12 illustrates the evolution of performance, measured in terms of MSE, of BM3D, DDID and their re-fitted versions as a function of the smoothing parameter $\gamma$. BM3D is studied on the *Monarch* image and DDID on the *Lena* image. Similar conclusions can be made from this experiment, re-fitting reaches its optimal performance for a larger smoothing parameter. Because the original estimator is nearly unbiased, re-fitting becomes challenging only when the original solution is over-smoothed, otherwise the gain in terms of bias is too small to compensate the loss in terms of variance. Given this high smoothing strength, some tiny structures have been lost, and thus the optimal re-fitting does not reach as good performance as the optimal original solution.

However, even tough the MSE is not necessarily improved, our re-fitting solutions present less artifacts (as confirmed by the very small gap in terms of SSIM values), see for instance the stripes of the *Monarch* or the left cheek of *Lena*. In fact, in order to recover details with low signal to noise ratio, the optimal original estimators authorize the apparition of low contrasted oscillating features. Nevertheless, a few of these oscillations tends to amplify some noise structures, hence, explaining these artifacts. In contrast, the optimal re-fitting prefers loosing such details rather than taking the risk of creating arbitrary structures and is thus more reliable.

Regarding our previous discussion, we believe that re-fitting would be nevertheless beneficial in terms of MSE for the class of images that are promoted by such restoration techniques. Characterizing this class of solutions for these methods is a very challenging topic out of the scope of this paper.

**6.6. Comparisons with other techniques devoted to the $\ell_1$ case.** We detail hereafter two different alternative strategies devoted to re-enhance the solution of the $\ell_1$ analysis regularization.

*Iterative hard-thresholding.* As shown earlier, the hard-thresholding is the re-fitted version of the soft-thresholding. Given an iterative solver $(k, y) \mapsto x^k$ composed of linear and soft-thresholding (such as primal-dual algorithms), one could consider replacing all soft-thresholding by hard-thresholding while keeping linearities unchanged: a technique often referred to as "iterative hard-thresholding" [2]. Unfortunately, such techniques only provide convergence to a local minimum of the $\ell_0$-regularized problem, and they do not converge to the sought re-fitting $\mathcal{R}_{\hat{x}}(y)$.

*Co-support identification based post re-fitting.* Another solution referred to as post re-fitting, and studied in, *e.g.,* [19, 36, 1, 28, 1], consists in identifying the (co-)support $\mathcal{I} = \{i \, : \, (\Gamma\hat{x}(y))_i \neq 0\}$ and solving a least-square problem constrained to $\{x \, : \, (\Gamma x)_{\mathcal{I}^c} = 0\}$, typically with conjugate gradient descent. However, $\hat{x}(y)$ is usually obtained thanks to a converging sequence $x^k$, and unfortunately, $\mathrm{supp}(\Gamma x^k)$ can be far from $\mathrm{supp}(\Gamma\hat{x}(y))$ even though $x^k$ can be made arbitrarily close to $\hat{x}(y)$. Such erroneous support identifications lead to results that strongly deviates from $\mathcal{R}_{\hat{x}}(y)$.

Figure 13 provides a comparison of our re-fitting method with the two other approaches mentioned earlier for aniso-TV used on an 8*bits* image (*peppers*), damaged by AWGN with $\sigma = 20$. The iterative hard-thresholding approach does not preserve the model space: transitions are not localized at the same positions as in the original solution and suspicious oscillations are created. The post re-fitting and our approach have both improved $\hat{x}(y)$ by enhancing each piece and preserving the location of transitions. Our method is nevertheless more stable than support identification which produces many errors due to wrong co-support identification.

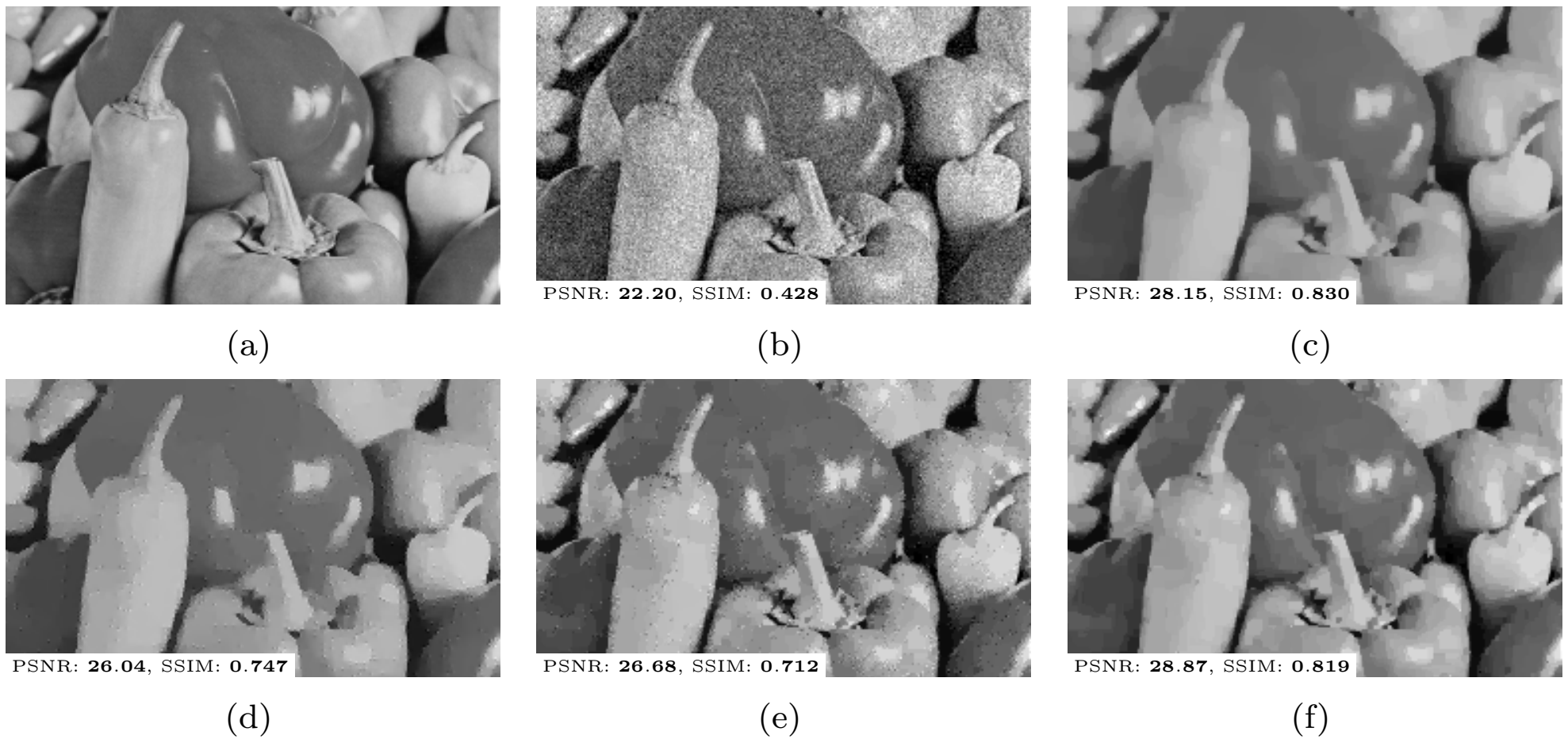

FIG. 13. *(a) Noise-free $x_0$. (b) Noisy $y = x_0 + w$. (c) Original aniso-TV $\hat{x}(y)$ with $\lambda = 1.2$. Enhanced results by (d) iterative hard-thresholding, (e) post re-fitting with support identification, and (f) our proposed covariant re-fitting.*

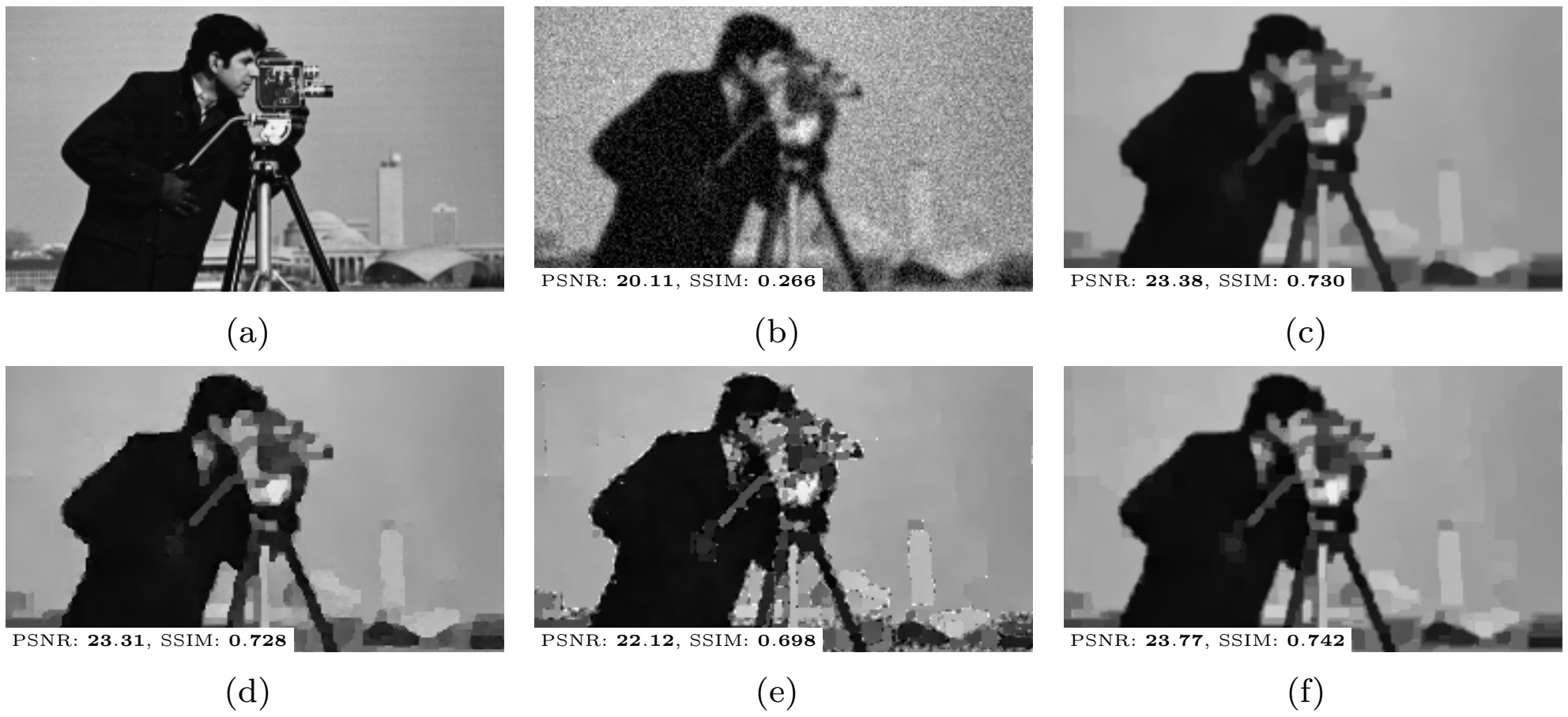

FIG. 14. *(a) Original $x_0$. (b) Blurred and noisy $y = \Phi x_0 + w$. (c) Original aniso-TV $\hat{x}(y)$ with $\lambda = 0.7$. Enhanced results by (d) iterative hard-thresholding, (e) post re-fitting with support identification, and (f) our proposed covariant re-fitting.*

Figure 14 provides another illustration highlighting the problem of support identification in an ill-posed problem. It consists of an $8bits$ image (*Cameraman*) damaged by a Gaussian blur of 2px and AWGN with $\sigma = 20$. Again, while aniso-TV reduces the contrast, the re-fitting recovers the original amplitudes and leaves discontinuities unchanged. Post re-fitting offers comparable results to ours except for suspicious oscillations due to wrong co-support identification.

In contrast with the support identification, CLEAR does not require the identification of the co-support, nor the identification of the model subspace $\mathcal{M}_{\hat{x}}(y)$. This is appealing since the co-support of $\hat{x}(y)$ is difficult to identify, particularly in the analy-sis context. Being computed along the iterations of the original algorithm, *i.e.,* jointly with $\hat{x}(y)$, our re-fitting strategy also provides more stable solutions.

**6.7. Comparisons with boosting strategies.** We detail hereafter other popular alternatives designed to re-enhance results of an arbitrary estimator.

**Twicing and boosting.** The boosting iterations introduced in [4] is a simple approach that consists in re-injecting to the current solution $\tilde{x}^k$ a filtered version of its residual $y - \Phi\tilde{x}^k$. The idea is that if parts of the signal where lost at iteration $k$, they might be retrieved in the residual. Given $\tilde{x}^0 = 0$, the iterations reads

$$\tilde{x}^{k+1} = \tilde{x}^k + \hat{x}(y - \Phi\tilde{x}^k) \ . \tag{51}$$

The first iterate is $\tilde{x}^1 = \hat{x}(y)$, and $\tilde{x}^2$ is known as the twicing estimate [47]. Such approaches are popular in non-parametric statistics, *e.g.,* for kernel smoothing [9].

When $k$ increases, its bias tends to decreases while its variance increases, see [43]. In denoising (*i.e.,* when $\Phi = \mathrm{Id}$) with a linear estimator $\hat{x}(y) = Wy$ (*e.g.,* the Tikhonov regularization), we get $\tilde{x}^k = (\mathrm{Id} - (\mathrm{Id} - W)^k)y$, for $k > 0$. In particular, the twicing reads as $\tilde{x}^2 = (2W - W^2)y$. In this linear case, the covariant re-fitting reads as $\mathcal{R}_{\hat{x}}(y) = (W + \rho W - \rho W^2)y$ and coincides with the twicing when $\rho = 1$.

**Iterative Bregman refinement.** In [33], the authors proposed an iterative procedure, originally designed to improve iso-TV results, given by

$$\tilde{x}^{k+1} = \hat{x}\left(y + \sum_{i=1}^{k}(y - \Phi\tilde{x}^k)\right) \ . \tag{52}$$

Unlike boosting that iteratively filters the residual, the idea is to filter a modified version of the input $y$ amplified by adding the sum of the residuals. When $\Phi = \mathrm{Id}$ and $\hat{x}(y) = Wy$, the iterative Bregman refinement reads as $\tilde{x}^k = (\mathrm{Id} - (\mathrm{Id} - W)^k)y$ and coincides with boosting (we refer to [5, 52, 23, 34] for related approaches).

**SOS-boosting.** In [37], the authors follow a similar idea by iteratively filtering a strengthened version of the input $y$. Their method, named Strengthen Operate Substract boosting (SOS-boosting), originally proposed for $\Phi = \mathrm{Id}$, performs iteratively the following update

$$\tilde{x}^{k+1} = \tau\hat{x}(y + \alpha\Phi\tilde{x}^k) - (\tau\alpha + \tau - 1)\tilde{x}^k \ . \tag{53}$$

where $\alpha$ and $\tau$ are two real parameters. The first one controls the emphasis of the solution (and the convergence of the sequence), while the second one controls the rate of convergence. When $\Phi = \mathrm{Id}$ and $\hat{x}(y) = Wy$, the SOS refinement with $\tau = 1$ reads as $\tilde{x}^k = Wy + \alpha(W - \mathrm{Id})\tilde{x}^{k-1}$, and in particular, for $k = 2$, we get $\tilde{x}^2 = (W - \alpha W + \alpha W^2)y$ which coincides with our covariant re-fitting for the choice $\alpha = -\rho$. For all considered estimators, we have always observed $\rho > 0$, contrarily to [37], where $\alpha > 0$ is implicitly assumed. Hence, we cannot concludes that the two models match in a specific linear setting. Another difference is that while we provide an automatic way to compute $\rho$ (see Equation (43)), the $\alpha$ parameter of the SOS-boosting must be tuned by the practitioner, a possibly cumbersome task, *e.g.,* when using cross-validation on a fix dataset of images and/or for varying noise levels.

**SAIF-boosting.** As described in [30], the diffusion of a filter consists in iteratively re-applying the filter to the current estimate $\tilde{x}^{k+1} = \hat{x}\left(\tilde{x}^k\right)$. The authors of [43] noticed that, unlike the boosting method of [4], the bias of this estimator increases and its variance decreases with $k$. As a consequence, the authors suggest mixing the two approaches by deciding at each iteration between performing a boosting or a diffusion step. To that end, they proposed a plug-in risk estimator that crudely estimates the MSE from a pre-filtered version of $y$. This approach is in fact applied locally on image patches, and is referred to as the Spatially Adapted Iterative Filter (SAIF)-boosting. Unlike the other techniques, the SAIF-boosting cannot be used as a black-box. Indeed, it requires to perform an eigen decomposition of $\hat{x}$ locally for each patch of $y$. This can be efficiently done for some kernel-based averaging filters, but can be very challenging for arbitrary estimators, such as for instance iso-TV.

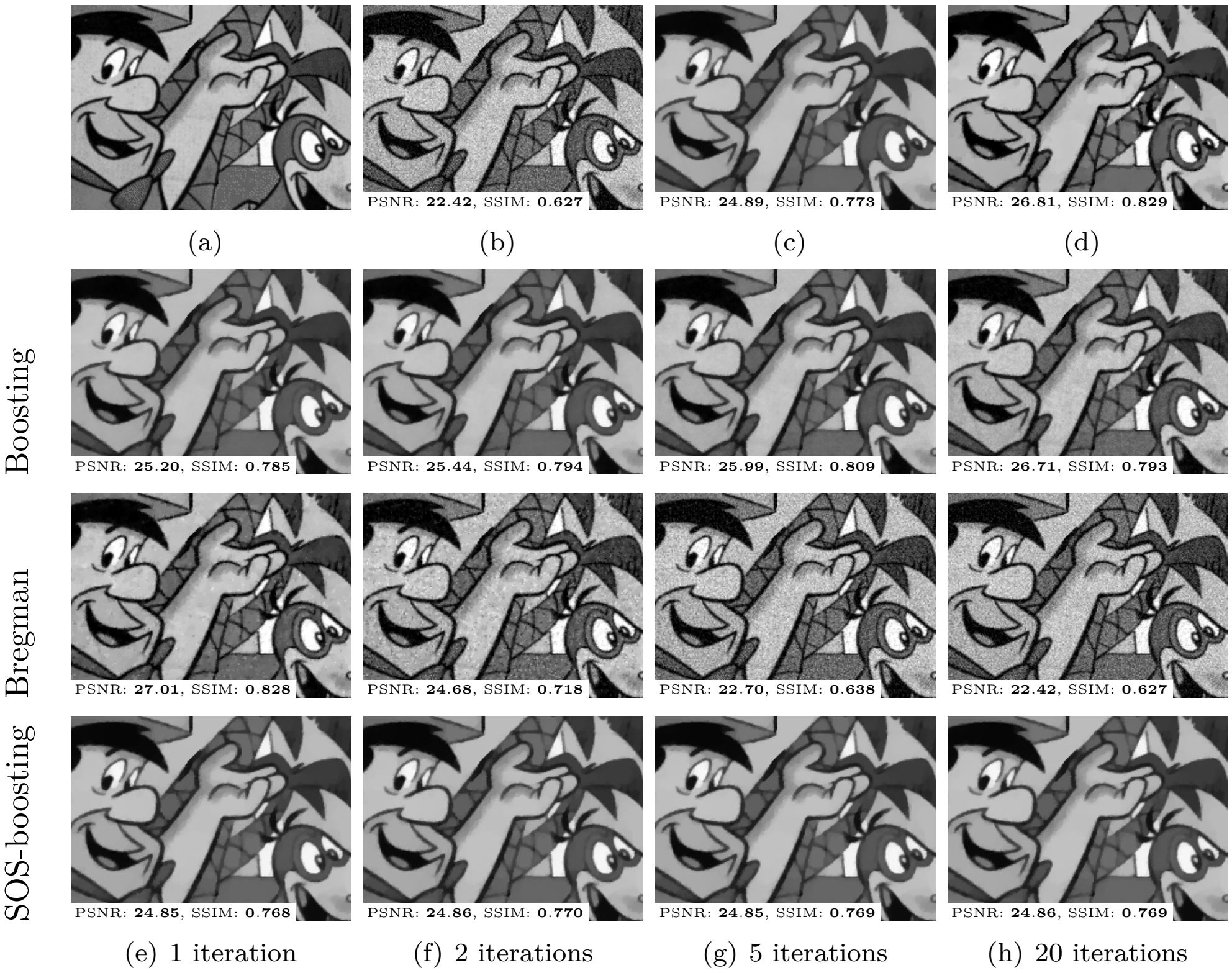

FIG. 15. *(a) Cartoon image* $x_0$. *(b) Noisy version* $y = x_0+w$. *(c) Iso-TV* $\hat{x}(y)$ *with* $\lambda = 1.5$. *(d) Our covariant re-fitting* $\mathcal{R}_{\hat{x}}(y)$. *(e,f,g,h) From top to bottom, boosting [47, 4], Bregman iterations [5] and SOS-boosting [37] at respectively* 1, 2, 5 *and* 20 *iterations.*

Though we have compared the expressions of boosting approaches with our re-fitting in the case of linear estimators, it is worth mentioning that boosting approaches are scarcely used in this case. Boosting appears more relevant in the non-linear case, since the successive re-application of a non-linear estimator $\hat{x}$ allows to recover parts of the signal that were lost at former iterations. Nevertheless, to boost the solution, the internal parameters of $\hat{x}$ may need to be re-adapted at each iteration, leading to cumbersome tuning of parameters in practice. Unlike boosting methods, a re-fitting approach should not modify the regularity and the structure of the first estimate. This is why CLEAR only considers the linearization of $\hat{x}$ at $y$ through the Jacobian.

It is important to have in mind that most of the theoretical results regarding boosting methods are often well grounded in the case where, even though $\hat{x}$ is non-linear, it acts locally as an averaging filter. In other words, locally, there exists a row stochastic linear operator $W$, *i.e.,* $W 1_n = 1_n$, (or even bi-stochastic, symmetric or independent of $y$) such that $\hat{x}(y) = Wy$. Such theoretical results could not be applied to the soft-thresholding (nor to more advanced methods we have considered). Indeed, for the soft-thresholding, a candidate for $W$ is the diagonal matrix defined as

$$W_{ii} = \begin{cases} 1 - \frac{\lambda}{|y_i|} & \text{if} \quad |y_i| > \lambda \ , \\ 0 & \text{otherwise} \ , \end{cases} \tag{54}$$

which is not row stochastic. In this context, the matrix $W$ is not the Jacobian, which is given in this case by the diagonal matrix

$$J_{ii} = \begin{cases} 1 & \text{if} \quad |y_i| > \lambda \ , \\ 0 & \text{otherwise} \ , \end{cases} \tag{55}$$

which, unlike $W$, is row stochastic and locally a projector. A second limitation is that even though $W$ is row stochastic, it might still encode a bias part. Typically, for the $\ell_1$ analysis described in Equation (5), a candidate for $W$ is

$$W = U(\Phi U)^+ - U(U^\top \Phi^\top \Phi U)^{-1} U^\top (\Gamma^\top)_{\mathcal{I}} R \ , \tag{56}$$

where $R$ is a diagonal matrix with diagonal elements $R_{ii} = \lambda/|y_i|$ for $i \in \mathcal{I}$ and 0 otherwise. One can check that if $1_n \in \mathrm{Ker}[\Gamma]$ (which holds true for aniso-TV), then $W$ is row stochastic. However, as seen in Equation (5), the quantity $U(U^\top \Phi^\top \Phi U)^{-1} U^\top (\Gamma^\top)_\mathcal{I} R$ is the term responsible for the systematic contraction of the $\ell_1$ analysis regularization (this simplifies to $\lambda/|y_i|$ for the soft-thresholding). As a consequence, the bias cannot be corrected by a single application of $W$. The Jacobian $J = U(\Phi U)^+$, which is again row stochastic, is free of this contraction term. Therefore our covariant re-fitting gets rid of this bias term after one single application of the Jacobian $J$.

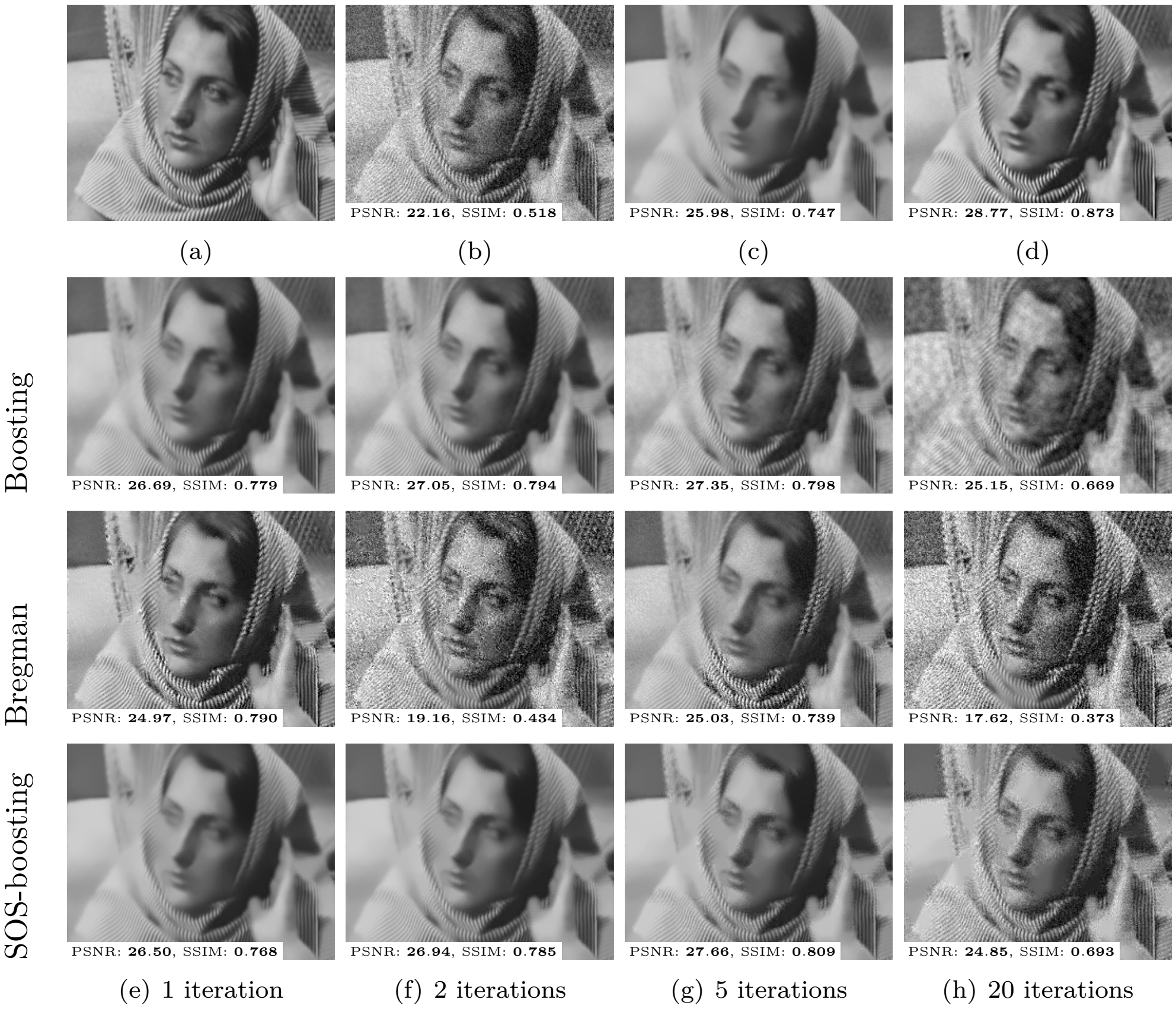

Fig. 16. *(a) Image $x_0$ with moderate patch redundancy. (b) Noisy version $y = x_0 + w$. (c) Non-local means $\hat{x}(y)$ with $h = 5.5$. (d) Covariant re-fitting $\mathcal{R}_{\hat{x}}(y)$. (e,f,g,h) From top to bottom, boosting [47, 4], Bregman iterations [5] and SOS-boosting [37] at respectively* 1, 2, 5 *and* 20 *iterations.*

Figure 15 and Figure 16 provide a comparison of our covariant re-fitting with the three first aforementioned boosting approaches, on two 8*bits* images (*Flinstones* and *Barbara*) damaged by AWGN with $\sigma = 20$, respectively for iso-TV and for non-local means. Note that the $\alpha$ and $\tau$ parameters of the SOS-boosting approach have been tuned to offer the most satisfying results, even though, we did not observe a significant impact in the iso-TV case. As expected, our covariant re-fitting provides results re-enhanced towards the amplitudes of the noisy inputs. In contrast, boosting approaches do not systematically re-fit towards the original amplitudes. While CLEAR preserves the structural content and smoothness of the original solution, the boosting approaches re-inject structural contents that were not originally preserved, and can present a large amount of residual noise. In these experiments, the original estimators were significantly biased, but for smaller parameter choices, re-injecting structural contents would have improved more the PSNR than re-fitting the amplitudes only. In fact, boosting techniques are more relevant to improve the quality of near unbiased estimates, while re-fitting techniques corrects only for the in-accuracy of biased but precise estimates (*i.e.,* with low estimation variance).

**7. Conclusion.** We have introduced a generalization of the popular least-square re-fitting, originally aimed at reducing the systematic contraction of the Lasso.

Together with this generalization, a generic implementation has been given for a wide class of ill-posed linear problem solvers. This implementation requires a computational overload of at most a factor three compared to the original solver; a factor that can even be reduced to two for most popular estimators in image processing.

While the classical re-fitting is inspired by the standard Lasso debiasing step (*i.e.,* least-square re-fitting on the estimated support), our generalization leverages the Jacobian of the estimator and does not rely on the notion of support, and its unstable identification. In particular, the proposed implementation only requires chain rules and differentiating the considered solver, and in practice it has also the benefit of increased stability compared to naive implementations.

For estimators such as Tikhonov regularization, total-variation non-local means, BM3D or DDID, numerical experiments have demonstrated the efficiency of the CLEAR technique in retrieving correct intensities while respecting the structure of the original biased estimator. Moreover, it has been shown in practice that re-fitting is beneficial when the underlying signal structure is well captured by the original estimator. Otherwise, re-fitting leads to too simplistic approximations, typically reflecting an inaccurate prior model. In other words, if the considered estimator is adequate with respect to the application context, then re-fitting is recommended.

We have highlighted the importance in distinguishing boosting approaches from the re-fitting one. In particular, re-fitting should be preferred in applications where the content of the original solution must be preserved. While boosting approaches are mostly used to enhance near unbiased estimators (typically coming from combinatorial or non-convex problems), the re-fitting is all the more relevant for estimators that present biases. For instance, re-fitting is essential for estimators solution of a convex problem that require a large bias correction to accurately retrieve the content of the signal of interest, a canonical example being the isotropic total-variation. Nevertheless, we believe that the notion of Jacobian based re-fitting could be of interest for boosting applications and we leave this to future work.

**7.1. Acknowledgment.** The authors would like to acknowledge the financial support of the GDR 720 ISIS (Information, Signal, Image et viSion). This study has also been carried out with financial support from the French State, managed by the French National Research Agency (ANR) in the frame of the "Investments for the future" Programme IdEx Bordeaux - CPU (ANR-10- IDEX-03-02).

**Appendix A. Sketch of proofs.** This section details how to retrieve closed-form expressions of some of the estimators studied in the paper.

**A.1. Retrieving the least-square solution.** We aim at retrieving here a minimizer of $\|\Phi x - y\|_2^2 + \iota_C(x)$ where $C = b + \mathrm{Im}[A]$, $b \in \mathbb{R}^p$ and $A \in \mathbb{R}^{p\times n}$. The initial problem can be recast as

$$\underset{x\in C}{\operatorname{argmin}}\ \|\Phi x - y\|_2^2 = b + A \cdot \underset{t\in\mathbb{R}^n}{\operatorname{argmin}}\ \|\Phi A t - (y - \Phi b)\|_2^2\ . \tag{57}$$

The right hand side problem being differentiable and convex, its first order optimality conditions give $A^\top \Phi^\top \Phi A t = A^\top \Phi^\top (y - \Phi b)$. In particular $t = (\Phi A)^+(y - \Phi b)$ is a solution, and hence $x = b + At$, *i.e.*, $x = b + A(\Phi A)^+(f - \Phi b)$ is a solution.

**A.2. Retrieving the Tikhonov solution.** We consider the minimization problem, defined for $\Gamma \in \mathbb{R}^{m\times p}$ and $\lambda > 0$, of $\frac{1}{2}\|\Phi x - y\|_2^2 + \frac{\lambda}{2}\|\Gamma x\|_2^2$. The objective being differentiable and convex, its first order optimality conditions give

$$\Phi^\top(\Phi x - y) + \lambda \Gamma^\top \Gamma x = 0 \Leftrightarrow (\Phi^\top \Phi + \lambda \Gamma^\top \Gamma) x = \Phi^\top f\ . \tag{58}$$

Provided $\mathrm{Ker}\,\Phi \cap \mathrm{Ker}\,\Gamma = \{0\}$ and $\lambda > 0$, the quantity $\Phi^\top\Phi + \lambda\Gamma^\top\Gamma$ is invertible and $x = (\Phi^\top\Phi + \lambda\Gamma^\top\Gamma)^{-1}\Phi^\top y$ is the unique solution.

**A.3. Retrieving the hard-thresholding solution.** We consider the minimization problem, defined for $\lambda > 0$, of $E(x, y) = \frac{1}{2}\|x - y\|_2^2 + \frac{\lambda^2}{2}\|x\|_0$. The problem is separable meaning that $\left[\underset{x\in\mathbb{R}^n}{\operatorname{argmin}}\ E(x,y)\right]_i = \underset{x_i\in\mathbb{R}}{\operatorname{argmin}}\ E_i(x_i, y_i)$ with

$$E_i(x_i, y_i) = \frac{1}{2}\begin{cases} y_i^2 & \text{if} \quad x_i = 0\ , \\ (x_i - y_i)^2 + \lambda^2 & \text{otherwise}\ . \end{cases} \tag{59}$$

Since $y_i^2 \leqslant \min_{x_i}[(x_i - y_i)^2 + \lambda^2] \Leftrightarrow |y_i| \leqslant \lambda$, we get

$$\min_{x_i} E_i(x_i, y_i) = \frac{1}{2}\left\{\begin{array}{ll} y_i^2 & \text{if} \quad |y_i| \leqslant \lambda\ , \\ \lambda^2 & \text{otherwise}\ , \end{array}\right. \tag{60}$$

which is reached by setting $x_i = 0$ when $|y_i| \leqslant \lambda$ and $x_i = y_i$ otherwise.

**A.4. Retrieving the soft-thresholding solution.** We consider the minimization problem, defined for $\lambda > 0$, of $E(x, y) = \frac{1}{2}\|x - y\|_2^2 + \lambda\|x\|_1$, which is, as for the hard-thresholding, separable with $E_i(x_i, y_i) = \frac{1}{2}(x_i - y_i)^2 + \lambda|x_i|$. By convexity, a minimum is reached when zero belongs to its sub-differential. Hence, $x_i$ is solution if

$$0 \in \partial E_i(x_i, y_i) \Leftrightarrow x_i \in y_i - \lambda \left\{\begin{array}{ll} \mathrm{sign}(x_i) & \text{if} \quad |x_i| > 0\ , \\ [-1, 1] & \text{otherwise}\ , \end{array}\right. \tag{61}$$

which holds true by setting $x_i = 0$ when $|y_i| \leqslant \lambda$ and $x_i = y_i - \lambda\,\mathrm{sign}(y_i)$ otherwise.

**A.5. Retrieving the (block-wise) non-local means.** We assume periodical boundary conditions such that all quantities $q$ indexed by $i = (i_1, i_2) \in \mathbb{Z}^2$ satisfies $q_i = q_{(i_1,i_2)} = q_{(i_1+k_1p_1, i_2+k_2p_2)}$ for all $(k_1, k_2) \in \mathbb{Z}^2$. This leads, for $l \in [-b, b]^2$, to the following relationship

$$
\begin{aligned}
F(x,y) &= \frac{1}{2}\sum_{i,j} w_{i,j}\|\mathcal{P}_i x - \mathcal{P}_j y\|_2^2 = \frac{1}{2}\sum_{i,j} w_{i,j}\sum_l (x_{i+l} - y_{j+l})^2 \\
&= \frac{1}{2}\sum_{i,j}\left[\sum_l w_{i,j}(x_{i+l} - y_{j+l})^2\right] = \frac{1}{2}\sum_{i,j}\left[\sum_l w_{i-l,j-l}(x_i - y_j)^2\right] \\
&= \frac{1}{2}\sum_{i,j}\left[\sum_l w_{i-l,j-l}\right](x_i - y_j)^2 = \frac{1}{2}\sum_{i,j}\bar{w}_{i,j}(x_i - y_j)^2 \ .
\end{aligned} \tag{62}
$$

For all $i \in [p_1] \times [p_2]$, studying the first optimality conditions gives

$$
\frac{\partial F(x,y)}{\partial x_i} = 0 \Leftrightarrow \sum_j \bar{w}_{i,j}(x_i - y_j) = 0 \Leftrightarrow x_i = \frac{\sum_j \bar{w}_{i,j} y_j}{\sum_j \bar{w}_{i,j}} \ . \tag{63}
$$

## Appendix B. Proof of Theorem 20.

Before turning to the proof of this theorem, let us introduce a first lemma.

LEMMA 22. *For all $y$, let $\hat{x}(y)$ be a solution of*

$$
\hat{x}(y) \in \operatorname*{argmin}_x F(y - \Phi x) + G(x) \ , \tag{64}
$$

*with $F$, $G$ two convex functions and $G$ being 1-homogeneous. Then for all $\varepsilon \in [0, 1]$, the following holds*

$$
(1-\varepsilon)\hat{x}(y) \in \operatorname*{argmin}_x F(y - \varepsilon\Phi\hat{x}(y) - \Phi x) + G(x) \ . \tag{65}
$$

*Proof.* Note that if $G$ is a convex and 1−homogeneous function, then $G$ is sub-additive, *i.e.*, $G(a)+G(b) \geqslant G(a+b)$. Next, assume that, for some $\varepsilon \in [0, 1]$, Eq. (65) does not hold, so that there exists $v$ such that

$$
\begin{aligned}
F(y - \varepsilon\Phi\hat{x}(y) - \Phi v) + G(v) &< F(y - \varepsilon\Phi\hat{x}(y) - (1-\varepsilon)\Phi\hat{x}(y)) + G((1-\varepsilon)\hat{x}(y)) \ , \\
&< F(y - \Phi\hat{x}(y)) + (1-\varepsilon)G(\hat{x}(y)) \ .
\end{aligned} \tag{66}
$$

It follows that

$$
F(y - \varepsilon\Phi\hat{x}(y) - \Phi v) + G(v) + \varepsilon G(\hat{x}(y)) < F(y - \hat{x}(y)) + G(\hat{x}(y)) \ . \tag{67}
$$

We also have $G(v) + \varepsilon G(\hat{x}(y)) = G(v) + G(\varepsilon\hat{x}(y)) \geqslant G(v + \varepsilon\hat{x}(y))$, since $G$ is 1-homogeneous and sub-additive. Hence, for $w = v + \varepsilon\hat{x}(y)$, we get

$$
F(y - \Phi w) + G(w) < F(y - \Phi\hat{x}(y)) + G(\hat{x}(y)) \ , \tag{68}
$$

which contradicts $\hat{x}(y) \in \operatorname{argmin} F(y - \Phi x) + G(x)$, and then concludes the proof. □

*Proof of Theorem 20.* By virtue of Lemma 22 and definition of $\hat{x}(y)$, we have $(1-\varepsilon)\hat{x}(y) = \hat{x}(y - \varepsilon\Phi\hat{x}(y))$ since $\hat{x}(y)$ is supposed to be the unique solution for all $y$. Now, recall that the linear Jacobian operator applied to $\Phi\hat{x}(y)$ is the directional derivative of $\hat{x}(y)$ in the direction $\Phi\hat{x}(y)$, then, for almost all $y$, we get

$$
J\Phi\hat{x}(y) \triangleq \lim_{\varepsilon\to 0}\frac{\hat{x}(y) - \hat{x}(y - \varepsilon\Phi\hat{x}(y))}{\varepsilon} = \lim_{\varepsilon\to 0}\frac{\hat{x}(y) - (1-\varepsilon)\hat{x}(y)}{\varepsilon} = \hat{x}(y) \ , \tag{69}
$$

which concludes the proof. □

**Appendix C. Proof of Theorem 21.**

Before turning to the proof of this theorem, let us introduce a first lemma.

Lemma 23. *The re-fitting $\mathcal{R}_{\hat{x}}(y)$ of the $\ell_1$ analysis regularization $x^\star = \hat{x}(y)$ is the solution of the saddle-point problem*

$$\min_{\tilde{x} \in \mathbb{R}^p} \max_{\tilde{z} \in \mathbb{R}^m} \|\Phi \tilde{x} - y\|_2^2 + \langle \Gamma \tilde{x}, \tilde{z} \rangle - \iota_{S_\mathcal{I}}(\tilde{z}) \ , \tag{70}$$

*where $\iota_{S_\mathcal{I}}$ is the indicator function of the convex set $S_\mathcal{I} = \{p \in \mathbb{R}^m \ : \ p_\mathcal{I} = 0\}$.*

*Proof.* As $\Phi U$ has full column rank, the re-fitting of the solution (5) is the unique solution of the constrained least-square problem (see Subsection 3.2)

$$\mathcal{R}_{\hat{x}}(y) = U(\Phi U)^+ y = \underset{\tilde{x} \in \mathcal{M}_{\hat{x}}(y)}{\operatorname{argmin}} \ \|\Phi \tilde{x} - y\|_2^2 \ . \tag{71}$$

Remark that $\tilde{x} \in \mathcal{M}_{\hat{x}}(y) = \operatorname{Ker}[\operatorname{Id}^t_{\mathcal{I}^c}\Gamma] \Leftrightarrow (\Gamma\tilde{x})_{\mathcal{I}^c} = 0 \Leftrightarrow \iota_{S_{\mathcal{I}^c}}(\Gamma\tilde{x}) = 0$, where $S_{\mathcal{I}^c} = \{p \in \mathbb{R}^m \ : \ p_{\mathcal{I}^c} = 0\}$.

Using Fenchel transform, $\iota_{S_{\mathcal{I}^c}}(\Gamma\tilde{x}) = \max_{\tilde{z}} \langle \Gamma\tilde{x}, \tilde{z}\rangle - \iota^*_{S_{\mathcal{I}^c}}(\tilde{z})$, where $\iota^*_{S_{\mathcal{I}^c}}$ is the convex conjugate of $\iota_{S_{\mathcal{I}^c}}$. Observing that $\iota_{S_\mathcal{I}} = \iota^*_{S_{\mathcal{I}^c}}$ concludes the proof. ∎

Given Lemma 23, replacing $\Pi_{z^k + \sigma\Gamma v^k}$ in (47) by the projection onto $S_\mathcal{I}$, *i.e.*,

$$\Pi_{S_\mathcal{I}}(\tilde{z})_{\mathcal{I}^c} = \tilde{z}_{\mathcal{I}^c} \quad \text{and} \quad \Pi_{S_\mathcal{I}}(\tilde{z})_\mathcal{I} = 0 \ , \tag{72}$$

leads to the primal-dual algorithm of [7] applied to problem (70) which converges to the re-fitted estimator $\mathcal{R}_{\hat{x}}(y)$. It remains to prove that the projection $\Pi_{z^k+\sigma\Gamma v^k}$ defined in (47) converges to $\Pi_{S_\mathcal{I}}$ in finite time.

*Proof of Theorem 21.* First consider $i \in \mathcal{I}$, *i.e.*, $|\Gamma x^\star|_i > 0$. By assumption on $\alpha$, $|\Gamma x^\star|_i \geq \alpha > 0$. Necessary $z^\star_i = \lambda \operatorname{sign}(\Gamma x^\star)_i$ in order to maximize (44). Hence, $|z^\star + \sigma\Gamma x^\star|_i \geq \lambda + \sigma\alpha$. Using the triangle inequality shows that

$$\lambda + \sigma\alpha \leqslant |z^\star + \sigma\Gamma x^\star|_i \leqslant |z^\star - z^k|_i + \sigma|\Gamma x^\star - \Gamma v^k|_i + |z^k + \sigma\Gamma v^k|_i \ . \tag{73}$$

Choose $\varepsilon > 0$ sufficiently small such that $\sigma\alpha - \varepsilon(1+\sigma) > \beta$. From the convergence of the primal-dual algorithm of [7], the sequence $(z^k, x^k, v^k)$ converges to $(z^\star, x^\star, x^\star)$. Therefore, for $k$ large enough, $|z^\star - z^k|_i < \varepsilon$, $|\Gamma x^\star - \Gamma v^k|_i < \varepsilon$, and

$$|z^k + \sigma\Gamma v^k|_i \geqslant \lambda + \sigma\alpha - \varepsilon(1+\sigma) > \lambda + \beta \ . \tag{74}$$

Next consider $i \in \mathcal{I}^c$, *i.e.*, $|\Gamma x^\star|_i = 0$, where by definition $|z^\star|_i \leqslant \lambda$. Using again the triangle inequality shows that

$$|z^k + \sigma\Gamma v^k|_i \leqslant |z^k - z^\star|_i + \sigma|\Gamma v^k - \Gamma x^\star|_i + |z^\star|_i \ . \tag{75}$$

Choose $\varepsilon > 0$ sufficiently small such that $\varepsilon(1+\sigma) < \beta$. As $(z^k, x^k, v^k) \to (z^\star, x^\star, x^\star)$, for $k$ large enough, $|z^k - z^\star|_i < \varepsilon$, $|\Gamma v^k - \Gamma x^\star|_i < \varepsilon$, and

$$|z^k + \sigma\Gamma v^k|_i < \lambda + \varepsilon(1+\sigma) \leqslant \lambda + \beta \ . \tag{76}$$

It follows that for $k$ sufficiently large $|z^k + \sigma\Gamma v^k|_i \leqslant \lambda + \beta$ if and only if $i \in \mathcal{I}^c$, and hence $\Pi_{z^k+\sigma K v^k}(\tilde{z}) = \Pi_{S_\mathcal{I}}(\tilde{z})$. As a result, all subsequent iterations of (47) will solve (70), and hence from Lemma 23 this concludes the proof of the theorem. ∎